\documentclass[journal]{IEEEtran}

\ifCLASSINFOpdf
\else
\fi
\usepackage{cite}
\usepackage{amssymb,amsfonts}
\usepackage{amsmath}
\usepackage{graphicx}
\usepackage{epstopdf}
\usepackage{color}
\usepackage{cases}
\usepackage{subfigure}
\usepackage{bm}

\usepackage{textcomp}
\usepackage{algorithm}
\usepackage{algorithmicx}

\newtheorem{assumption}{Assumption}
\newtheorem{lemma}{Lemma}
\newtheorem{theorem}{Theorem}

\newtheorem{remark}{Remark}

\newtheorem{corollary}{Corollary}

\begin{document}
%
\title{Decentralized Nash Equilibria Learning for Online Game with Bandit Feedback}
%
%
%

\author{
Min Meng, Xiuxian Li,  {\em Senior Member, IEEE}, and Jie Chen, {\em Fellow, IEEE}
\thanks{This work was partially supported by Shanghai Pujiang Program under Grand 21PJ1413100, the National Natural Science Foundation of China under Grand 62003243 and 62103305, Shanghai Municipal Science and Technology Major Project under Grand 2021SHZDZX0100, Shanghai Municipal Commission of Science and Technology under Grand 19511132101, Young Elite Scientist Sponsorship Program by cast of China Association for Science and Technology under Grand YESS20200136, and the Basic Science Centre Program by the National Natural Science Foundation of China under Grant 62088101. (\emph{Corresponding author: Xiuxian Li.})}
\thanks{M. Meng, X. Li and J. Chen are with the Department of Control Science and Engineering, College of Electronics and Information Engineering, and Shanghai Research Institute for Intelligent Autonomous Systems, Tongji University, Shanghai, China (Email: mengmin@tongji.edu.cn; xli@tongji.edu.cn; chenjie206@tongji.edu.cn).
}
}

\maketitle

\begin{abstract}
This paper studies distributed online bandit learning of generalized Nash equilibria for online game, where cost functions of all players and coupled constraints are time-varying. The values rather than full information of cost and local constraint functions are revealed to local players gradually. The goal of each player is to selfishly minimize its own cost function with no future information subject to a strategy set constraint and time-varying coupled inequality constraints. To this end, a distributed online algorithm based on mirror descent and one-point bandit feedback is designed for seeking generalized Nash equilibria of the online game. It is shown that the devised online algorithm achieves sublinear expected regrets and accumulated constraint violation if the path variation of the generalized Nash equilibrium sequence is sublinear. Furthermore, the proposed algorithm is extended to the scenario of delayed bandit feedback, that is, the values of cost and constraint functions are disclosed to local players with time delays. It is also demonstrated that the online algorithm with delayed bandit feedback still has sublinear expected regrets and accumulated constraint violation under some conditions on the path variation and delay. Simulations are presented to illustrate the efficiency of theoretical results.
\end{abstract}


%
\IEEEpeerreviewmaketitle

\section{Introduction}
Multi-player noncooperative games with self-interested decision makers have found a remarkable breadth of applications, such as social networks \cite{ghaderi2014opinion}, smart grid \cite{saad2012game}, sensor networks \cite{stankovic2012distributed}, and so on. A vital concept for this kind of games is \emph{Nash Equilibrium} (NE) \cite{basar1999dynamic}, from which no player has an incentive to deviate. If the strategy set of each player depends on other players' strategies, which often emerges in a wide range of real-world applications, e.g., limited resource among all players, then the NE is called a generalized NE (GNE).

Recently, distributed NEs and GNEs seeking in noncooperative games have received increasing attention. In contrast to centralized methods \cite{facchinei2010generalized,shamma2005dynamic}, distributed NEs or GNEs seeking algorithms \cite{de2019distributed,gadjov2019passivity,koshal2016distributed,salehisadaghiani2019distributed,ye2021distributed}, without the need of a global coordinator bidirectionally communicating with all the players, only depend on partial players' decision instead of full actions' information. As a result, distributed algorithms can reduce communication burden, increase robustness to link failures or malicious attacks, and preserve individual players' private information to some extent.

%

By now the above discussions have been on offline game, where cost and constraint functions are time-invariant. Nevertheless, the surrounding environments in various practical situations, such as real-time traffic networks, online auction and allocation radio resources, often change over time, incurring time-varying cost functions and/or constraints, which is usually called online game. In online game, the cost and constraint functions are revealed to local players only after making their decisions. In this setting, a distributed GNE learning algorithm for online game with time-invariant constraints was designed in \cite{lu2020online} based on primal-dual strategies and distributed consensus. Subsequently, the authors of \cite{meng2021decentralized} considered time-varying constraints and proposed a distributed online algorithm on the basis of mirror descent and primal-dual strategies.

In \cite{lu2020online,meng2021decentralized}, each player can access the gradient information of its cost and local constraint functions sequentially. However, the gradient information in many realistic applications cannot be grabbed by local players, especially if the cost and constraint functions are not revealed. Instead, only function values, even with time-delays due to latency in communication and computation, are available to local players, i.e., bandit feedback, making online game more considerably challenging \cite{bravo2018bandit,mertikopoulos2019learning}.

In this paper, distributed online learning for GNE of online game with time-varying coupled constraints is investigated. Both delay-free and delayed bandit feedbacks are investigated. The main contributions of this paper are summarized as follows.
\begin{itemize}
\item[1)] To the best of our knowledge, this paper is the first to study online game with time-varying coupled inequality constraints with bandit feedback, while existing works on online game only consider the case with full information of cost and constraint functions \cite{lu2020online,meng2021decentralized}.
\item[2)] A distributed GNE seeking algorithm for online game is devised by mirror descent and one-point bandit feedback. It is demonstrated that this algorithm can achieve sublinear expected regrets and accumulated constraint violation if the path variation of the GNE sequence is sublinear.
\item[3)] The above algorithm is extended to the case with delayed bandit feedback. Also, it is shown that the generated expected regrets and accumulated constraint violation by the new algorithm are sublinear under some conditions on the path variation and delay.
\end{itemize}

The rest of this paper is organized as follows. In Section \ref{section2}, the problem formulation and some preliminaries on Bregman divergence, as well as one-point sampling gradient estimator, are introduced. Sections \ref{section3} and \ref{section4} propose distributed bandit online algorithms for seeking GNEs without and with delays, respectively. Numerical examples are presented to support the obtained results in Section \ref{section5}. Section \ref{section6} concludes this paper.

{\em Notations.} $\mathbb{R}$, $\mathbb{R}^n$ and $\mathbb{R}^{m\times n}$ represents the sets of real numbers, real column vectors of dimension $n$ and real matrices of dimension $m\times n$, respectively. $\mathbb{B}^n$ and $\mathbb{S}^n$ represents the unit ball and the sphere centered at the origin in $\mathbb{R}^n$, respectively. For a positive integer $m$, $[m]:=\{1,2,\ldots,m\}$. Denote by ${\bf 1}_n$ (resp. ${\bf 0}_n$) the $n$-dimensional vector with all elements being 1 (resp. 0). The transpose of a vector or matrix $P$ is denoted as $P^\top$. $col(z_1,\ldots,z_n):=(z_1^{\top},\ldots,z_n^{\top})^{\top}$. $A\otimes B$ is the Kronecker product of matrices $A$ and $B$. $\langle x,y\rangle$ represents the inner product of vectors $x$ and $y$. For two vectors/matrices $w,v$ of the same dimension, $w\leq v$ (resp. $w\geq v$) means that each entry of $w-v$ is nonpositive (resp. nonnegative), while for two real symmetric matrices $P$ and $Q$, $P\succeq~Q$ represents that $P-Q$ is positive semi-definite.  For $x\in\mathbb{R}^n$, $[x]_+$ denotes the projection of $x$ onto $\mathbb{R}^n_+$, where $\mathbb{R}_+^n:=\{y\in\mathbb{R}^n\mid y\geq0\}$. The symbol $h_1=\mathcal{O}(h_2)$ for two functions $h_1(\cdot)$ and $h_2(\cdot)$ means that there exists a positive constant $C$ such that $|h_1(x)|\leq Ch_2(x)$ and $h_1={\bf o}(T)$ connotes that $\lim\limits_{T\to\infty}\frac{h_1}{T}=0$.


\section{Problem Formulation and Preliminaries}\label{section2}
In this section, the studied problem is first introduced and then some necessary preliminaries on Bregman divergence and one-point sampling gradient estimator are presented.
\subsection{Problem Formulation}

The online game with time-varying coupled inequality constraints studied in this paper is denoted as $\Gamma(\mathcal{V},X_t,f_t)$. $\mathcal{V}:=[N]$ is the set of $N$ players. $f_{t}:=(f_{1,t},\ldots,f_{N,t})$, where $f_{i,t}:X_i\to\mathbb{R}$ is the private time-varying cost function of player $i$ at time step $t\geq1$. $X_t=X_{0,t}\bigcap(X_{1}\times\cdots\times X_{N})$ is the time-varying constraint, where $X_{0,t}:=\{x\in\mathbb{R}^n\mid g_t(x):=\sum_{i=1}^Ng_{i,t}(x_i)\leq{\bf 0}_m\}$ is time-varying coupled inequality constraints and $X_i\subseteq\mathbb{R}^{n_i}$ imposes the time-invariant private strategy set of player $i$. Here, $n:=\sum_{i=1}^Nn_i$ and $g_{i,t}:\mathbb{R}^{n_i}\to\mathbb{R}^m$ is the private constraint function of player $i$.

In online game, $f_{i,t}$ and $g_{i,t}$ are revealed to agent $i$ only after making its decision at time $t$, without any knowledge of future information. Full information on cost and constraint functions may be impossible or computationally heavy to be disclosed to local players in many applications. Hence, this paper focuses on the case where only function values of $f_{i,t}(x_{i,t},x_{-i,t})$ and $g_{i,t}(x_{i,t})$ are revealed privately to player $i$ after a strategy profile $x_t:=(x_{1,t},\ldots,x_{N,t})\in{X}_t$ is chosen based on historical information at hand, where $x_{i,t}\in{X}_i$ is the strategy selected by player $i$ at time $t$ and $x_{-i,t}:= col(x_{1,t}, \ldots, x_{i-1,t}, x_{i+1,t}, \ldots, x_{N,t})$ denotes the joint strategy of all players except player $i$ at time $t$. Also, each player is unable to acquire all other players' information, and instead has the ability to receive the information of partial players via local communications.

The communication pattern among all players is captured by a graph ${\mathcal G}=({\mathcal V},{\mathcal E},A)$, where ${\mathcal V}=[N]$ is the node set, ${\mathcal E}\subseteq{\mathcal V}\times{\mathcal V}$ is the edge set, and $A=(a_{ij})\in\mathbb{R}^{N\times N}$ is the weighted adjacency matrix. $a_{ij}>0$ if $(j,i)\in{\mathcal E}$ and $a_{ij}=0$ otherwise. $a_{ii}$ for all $i\in[N]$ are assumed to be positive  in this paper. The graph is undirected if $(i,j)\in\mathcal{E}$ implies $(j,i)\in\mathcal{E}$. $j$ in an edge $(j,i)$ is called a neighbor of $i$. Denote by ${\mathcal N}_i=\{j:~(j,i)\in{\mathcal E}\}$ the set of neighbors of node $i$. A path from node $i_1$ to node $i_l$ is a sequence of edges $(i_h,i_{h+1})$, $h=1,2,\ldots,l-1$, which joins a sequence of distinct vertices. An undirected graph ${\mathcal G}$ is said to be connected if there is a path between any two nodes.

For online game $\Gamma(\mathcal{V},X_t,f_t)$, it is impossible for agents to pre-compute a GNE $(x^*_{i,t},x^*_{-i,t})$, i.e., a strategy profile $(x^*_{i,t},x^*_{-i,t})\in X_t$ satisfying
\begin{align}\label{equ1}
f_{i,t}(x_{i,t}^*,x_{-i,t}^*)\leq f_i(x_{i,t},x_{-i,t}^*), ~\forall x_{i,t}\in X_{i,t}(x^*_{-i,t}),
\end{align}
 where $X_{i,t}(x_{-i,t}):=\{x_i\in X_i\mid(x_i,x_{-i,t})\in{X}_{t}\}$,
but agents can improve their performance by learning through play. Note that (\ref{equ1}) is equivalent to the following optimization problem:
\begin{align}\label{equ2}
&\min\limits_{x_{i}}~~~~~~ f_{i,t}(x_{i},x^*_{-i,t})\nonumber\\
&\text{subject to} ~~x_{i}\in X_{i,t}(x^*_{-i,t}).
\end{align}
Then in this paper, based on (\ref{equ2}), an important performance measure called \emph{dynamic regret}, is adopted, defined as follows:
\begin{align}\label{equ3}
Reg_i(T):=\sum_{t=1}^T(f_{i,t}(x_{i,t},x^*_{-i,t})-f_{i,t}(x^*_{i,t},x^*_{-i,t})),
\end{align}
where $T$ is the total learning time horizon.

In the meantime, to measure the violation of constraints, the commonly used constraint violation measure is presented as
\begin{align}\label{equ5}
R_g(T):=\left\|\left[\sum\limits_{t=1}^Tg_t(x_t)\right]_+\right\|.
\end{align}

The metrics in (\ref{equ3}) and (\ref{equ5}) provide a meaningful method for quantifying the ability of an online algorithm to adapt to unknown and unpredictable environments. From this viewpoint, an online algorithm is said to be ``good'' or  ``no-regret'' if $Reg_i(T)$ and $R_g(T)$ are sublinear with respect to $T$, i.e.,
\begin{align}
Reg_i(T)={\bf o}(T),~i\in[N],~~R_g(T)&={\bf o}(T).
\end{align}

However, if the GNE sequence $\{x_t^*\}_{t=1}^T$ fluctuate drastically, then learning the exact GNE may be impossible, as studied in online optimization \cite{zinkevich2003online}. To tackle this issue, the \emph{path variation} (or \emph{path length}) of the GNE sequence, $\Theta_T^*$, is usually employed, defined as \cite{hall2015online,lu2020online}:
\begin{align}
\Theta_T^*:=\sum_{t=1}^T\|x^*_{t+1}-x^*_t\|.
\end{align}

To proceed, some standard assumptions on the studied game are made \cite{salehisadaghiani2016distributed,lu2020online}.
\begin{assumption}\label{assump1}~
\begin{itemize}
\item[1)] The sets $X_i$, $i\in[N]$, are nonempty, compact and convex. In addition, two positive constants $r_i,R_i$ exist such that
\begin{align}
r_i\mathbb{B}^{n_i}\subseteq X_i\subseteq R_i\mathbb{B}^{n_i},~\forall i\in[N].
\end{align}
Here, $r_i$, $i\in[N]$, are known a priori.
\item[2)] For any given $x_{-i}\in\mathbb{R}^{n-n_i}$, $f_{i,t}(x_i,x_{-i})$ is differentiable and convex with respect to $x_i$, and $g_{ij,t}(x_i)$ is differentiable and convex, where $g_{ij,t}$ is the $j$th element of $g_{i,t}$, $j\in[m]$.
\item[3)] The constraint set $X_t$ is nonempty and also satisfies Slater's constraint qualification.
\end{itemize}
\end{assumption}

By Assumption \ref{assump1}, one can further assume the uniform boundedness, i.e., for any $t$, $i\in[N]$, $x_i\in X_{i}$ and $x_{-i}\in\mathbb{R}^{n-n_i}$, there exist positive constants $B_x$, $B_f$, $L_f$, $B_g$ and $L_g$ such that
\begin{align}
&\|x_i\|\leq B_x,\|f_{i,t}(x_i,x_{-i})\|\leq B_f,\|g_{i,t}(x_i)\|\leq B_g,\label{equ7}\\
&\|\nabla_if_{i,t}(x)\|\leq L_f,\|\nabla g_{i,t}(x_i)\|\leq L_g,\label{equ8}
\end{align} where $\nabla g_{i,t}(x_i):=(\nabla g_{i1,t}(x_i),\ldots,\nabla g_{im,t}(x_i))^{\top}\in\mathbb{R}^{m\times{n}_i}$ and $\nabla_if_{i,t}(x_i,x_{-i}):=\frac{\partial f_{i,t}(x_i,x_{-i})}{\partial x_i}$.

\begin{assumption}\label{assump2}
For any $i\in[N]$, $\nabla_if_{i,t}(x_{i},x_{-i})$ is $L$-Lipschitz continuous, i.e.,
\begin{align}
\|\nabla_if_{i,t}(x_{i},x_{-i})-\nabla_if_{i,t}(y_{i},y_{-i})\|\leq L\|x-y\|
\end{align} for any $x_i,y_i\in X_i$ and $x_{-i},y_{-i}\in\mathbb{R}^{n-n_i}$.
\end{assumption}

Define a pseudo-gradient mapping of $\Gamma(\mathcal{V},X_t,f_t)$ as
\begin{align}
F_t(x):=col(\nabla_1f_{1,t}(x),\ldots,\nabla_Nf_{N,t}(x)).
\end{align}

\begin{assumption}\label{assump3}
The mapping $F_t(x)$ is $\mu$-strongly monotone on the set $X_{1}\times\cdots\times X_{N}$, i.e., for any $x,y\in X_{1}\times\cdots\times X_{N}$,
\begin{align}
(F_t(x)-F_t(y))^{\top}(x-y)\geq\mu\|x-y\|^2,
\end{align}
which is equivalent to
\begin{align}
\sum_{i=1}^N(\nabla_if_{i,t}(x)-\nabla_if_{i,t}(y))^{\top}(x_i-y_i)\geq\mu\|x-y\|^2,
\end{align}
where $x=col(x_1,\ldots,x_N)$ and $y=col(y_1,\ldots,y_N)$.
\end{assumption}

Under Assumption \ref{assump1}, it is obtained from Theorem 3.9 in \cite{facchinei2009nash} that at any $t$, a solution $x_t^*\in X_t$ to the variational inequality
\begin{align}\label{equ14}
(F_t(x^*_t))^{\top}(x-x^*_t)\geq 0 ~\text{for all}~x\in X_t
\end{align}
is also a GNE of game $\Gamma(\mathcal{V},X_t,f_t)$, and such GNE $x^*_t$ is also called a variational GNE. Additionally, the solution of the variational inequality (\ref{equ14}) is unique under Assumption \ref{assump3}. It should be noted that seeking all GNEs is rather difficult even for offline game, and thereby this paper focuses on seeking the unique variational GNE sequence. As a matter of fact, seeking the variational GNE of games with coupled constraints was widely studied \cite{pavel2020distributed,lu2020online,meng2021decentralized} since the variational GNE has good stability with the economic interpretation of no price discrimination.

Some assumptions on players' communication are listed below.
\begin{assumption}\label{assump4}~
\begin{itemize}
\item[1)] The interaction graph ${\cal G}=([N],{\mathcal E},A)$ is undirected and connected.
\item[2)] $A^{\top}=A$ and $A{\bf1}_N={\bf1}_N$.
\end{itemize}
\end{assumption}

Define $$\sigma_m:=\left\|A-\frac{1}{N}{\bf1}_N{\bf1}_N^{\top}\right\|.$$ Under Assumption \ref{assump5}, it is known from \cite{li2020} that $\sigma_m\in[0,1)$.

\subsection{Bregman Divergence}
Each player $i\in[N]$ is assigned a differentiable function $\phi_i: X_i\to\mathbb{R}$, which is assumed to be $\mu_i$-strongly convex for some $\mu_i>0$, i.e.,
\begin{align*}
\phi_i(\theta)\geq\phi_i(\vartheta)+\langle\nabla\phi_i(\vartheta),\theta-\vartheta\rangle
+\frac{\mu_i}{2}\|\theta-\vartheta\|^2, \forall \theta,\vartheta\in X_i.
\end{align*}

The Bregman divergence $D_{\phi_i}(\theta,\vartheta)$ for two points $\theta,\vartheta\in{X}_i$ corresponding to $\phi_i$ is defined as
\begin{align}
D_{\phi_i}(\theta,\vartheta)
:=\phi_i(\theta)-\phi_i(\vartheta)-\langle\nabla\phi_i(\vartheta),\theta-\vartheta\rangle
\end{align}
with the properties below:
\begin{itemize}
\item[i)] $D_{\phi_i}(\cdot,\cdot)$ is $\mu_0$-strongly convex with respect to the first variable, where $\mu_0:=\min\{\mu_1,\ldots,\mu_N\}$.
\item[ii)] It holds
\begin{align}\label{equ17}
D_{\phi_i}(\theta,\vartheta)\geq\frac{\mu_0}{2}\|\theta-\vartheta\|^2,~\forall \theta,\vartheta\in X_i.
\end{align}
\item[iii)] The generalized triangle inequality is satisfied, i.e.,
\begin{align}
&\langle \theta-\vartheta, \nabla\phi_i(\vartheta)-\nabla\phi_i(\xi)\rangle\nonumber\\
&=D_{\phi_i}(\theta,\xi)-D_{\phi_i}(\theta,\vartheta)-D_{\phi_i}(\vartheta,\xi)\label{equ18}
\end{align}
for any $\theta,\vartheta,\xi\in X_i$.
\end{itemize}


Two mild assumptions on Bregman divergence are made.
\begin{assumption}\label{assump5}
$D_{\phi_i}(\cdot,\cdot)$ is Lipschitz with respect to the first variable, i.e., there is a constant $K>0$ such that for any $\theta_1,\theta_2,\vartheta\in X_{i}$,
\begin{align}
|D_{\phi_i}(\theta_1,\vartheta)-D_{\phi_i}(\theta_2,\vartheta)|\leq K\|\theta_1-\theta_2\|.
\end{align}
\end{assumption}
\begin{assumption}\label{assump6}
$D_{\phi_i}(\cdot,\cdot)$ is convex with the second variable, i.e., for any $\alpha\in[0,1]$ and $\vartheta_1,\vartheta_2,\theta\in X_{i}$,
\begin{align}
&D_{\phi_i}(\theta,\alpha\vartheta_1+(1-\alpha)\vartheta_2)\nonumber\\
&\leq\alpha D_{\phi_i}(\theta,\vartheta_1)+(1-\alpha)D_{\phi_i}(\theta,\vartheta_2).
\end{align}
\end{assumption}

Assumption \ref{assump5} is satisfied if $\phi_i(\theta)$ is Lipschitz on $X_{i}$. It can also be obtained from Assumption \ref{assump5} that for any $\theta,\vartheta\in X_i$,
\begin{align}\label{equ20}
D_{\phi_i}(\theta,\vartheta)=|D_{\phi_i}(\theta,\vartheta)-D_{\phi_i}(\vartheta,\vartheta)|\leq K\|\theta-\vartheta\|.
\end{align}
Assumption \ref{assump6} plays an important role in analyzing the designed algorithms and can be ensured when $\phi_i(\theta)$ is thrice continuously differentiable and
$H_{\phi_i}(\theta)\succeq 0$, $H_{\phi_{i}}(\theta)+\nabla H_{\phi_i}(\theta)(\theta-\vartheta)\succeq0$ for $\theta,\vartheta\in X_{i}$, where $H_{\phi_i}$ denotes the Hessian matrix of $\phi_i$ \cite{bauschke2001joint}.

\subsection{Gradient Approximation}
To facilitate the following algorithm development and performance analysis, this subsection briefly introduces a gradient approximation approach. Let $f:\mathbb{X}\to\mathbb{R}$ be a function, where the domain $\mathbb{X}\subseteq\mathbb{R}^n$ is a convex and bounded set, and has a nonempty interior. Assume that $\mathbb{X}$ is contained in the ball $R(\mathbb{X})\mathbb{B}^n$ and contains the ball $r(\mathbb{X})\mathbb{B}^n$, that is, $r(\mathbb{X})\mathbb{B}^n\subseteq\mathbb{X}\subseteq R(\mathbb{X})\mathbb{B}^n$, where $r(\mathbb{X})$ and $R(\mathbb{X})$ are positive constants. A one-point sampling gradient estimator is proposed as \cite{flaxman2005online}
\begin{align}\label{equ21}
\hat{\nabla}f(x)=\frac{n}{\delta}f(x+\delta u)u,~\forall x\in(1-\eta)\mathbb{X},
\end{align} where $u\in\mathbb{S}^n$ is a uniformly distributed random vector, $\eta\in(0,1)$ is a shrinkage coefficient and $\delta\in(0,r(\mathbb{X})\eta]$ is an exploration parameter.
It should be noted that the perturbations in (\ref{equ21}) can be ensured to still remain in the set $\mathbb{X}$ by defining the estimator (\ref{equ21}) over the shrinking set $(1-\eta)\mathbb{X}$ instead of $\mathbb{X}$. Let us define the uniformly smoothed version of $f$ as
$$\hat{f}(x):={\bf E}_{v\in\mathbb{B}^n}[f(x+\delta v)].$$
\begin{lemma}[\cite{yi2020distributed}]\label{lemma1}
\begin{itemize}
\item[1)] It holds that $x+\delta u\in\mathbb{X}$ for any $x\in(1-\eta)\mathbb{X}$, $u\in\mathbb{S}^n$ and $\delta\in(0,r(\mathbb{X})\eta]$.
\item[2)] $\hat{f}(x)$ is differentiable on $(1-\eta)\mathbb{X}$ even when $f(x)$ is not, and there holds that
\begin{align}
\nabla\hat{f}(x)={\bf E}_{u\in\mathbb{S}^n}\left[\hat{\nabla}f(x)\right],~\forall x\in(1-\eta)\mathbb{X}.
\end{align}
\item[3)] $\hat{f}(x)$ is convex on $(1-\eta)\mathbb{X}$ if $f(x)$ is convex on $\mathbb{X}$, and $f(x)\leq\hat{f}(x)$ for any $x\in(1-\eta)\mathbb{X}$.
\item[4)] If $f(x)$ is $l_0$-Lipschitz on the set $\mathbb{X}$, then $\hat{f}$ and $\hat{\nabla}f$ are $l_0$-Lipschitz and $nl_0/\delta$-Lipschitz, respectively. Additionally, there holds that
    \begin{align}
    |\hat{f}(x)-f(x)|\leq\delta l_0,~\forall x\in(1-\eta)\mathbb{X}.
    \end{align}
\item[5)] If $|f(x)|\leq F_0$ for any $x\in\mathbb{X}$, then
\begin{align}
|\hat{f}(x)|&\leq F_0,\\
\|\hat{\nabla}f(x)\|&\leq\frac{nF_0}{\delta},~\forall x\in(1-\eta)\mathbb{X}.
\end{align}
\end{itemize}
\end{lemma}

In fact, replacing the original function $f$ by the function $\hat{f}$ is the critical idea of this gradient-free method since it is seen from 4) of Lemma \ref{lemma1} that the function $f(x)$ is extremely close to $\hat{f}(x)$ when $\delta$ is small. In addition, $\hat{\nabla}f$ can be regarded as an estimator of $\nabla\hat{f}$ from 2) of Lemma \ref{lemma1}.

\section{Distributed Bandit Feedback}\label{section3}
In this section, a distributed online algorithm for tracking the variational GNE sequence of the studied online game is proposed based on one-point bandit feedback method and mirror descent. We also analyze the expected regrets and constraint violation bounds for the proposed algorithm.


\begin{algorithm}[!ht]
\caption{Distributed Online Primal-Dual Mirror Descent with One-point Bandit Feedback}\label{alg1}
Each player $i$ maintains vector variables $x_{i,t}\in{X_i}$, $z_{i,t}\in(1-\eta_{t})X_i$, $\tilde{z}_{i,t}\in(1-\eta_{t})X_i$ and $\lambda_{i,t}\in\mathbb{R}_+^{m}$ at iteration $t$.

 {\bf Initialization:} For any $i\in[N]$, initialize $z_{i,1},\tilde{z}_{i,1}\in(1-\eta_{1})X_i$ arbitrarily, $x_{i,1}=z_{i,1}+\delta_1u_{i,1}$, and $\lambda_{i,1}={\bf0}_m$.

{\bf Iteration:} At $t\geq 1$, every player $i$ receives the function values $f_{i,t}(x_{i,t},x_{-i,t})$ and $g_{i,t}(x_{i,t})$ after the strategy profile $(x_{i,t},x_{-i,t})$ is made, and performs the following update:
\begin{subequations}
\begin{align}
\tilde{z}_{i,t+1}&=\arg\min\limits_{z\in(1-\eta_t)X_i}\{\alpha_t\langle{z},\hat{\nabla}_if_{i,t}({z}_{t})\rangle\nonumber\\
&~~~+\alpha_t\langle{z},(\hat{\nabla}{g}_{i,t}(z_{i,t}))^{\top}\tilde{\lambda}_{i,t}\rangle+D_{\phi_i}(z,z_{i,t})\},\label{equ31b}\\
z_{i,t+1}&=(1-\alpha_t)z_{i,t}+\alpha_t\tilde{z}_{i,t+1},\label{}\\
x_{i,t+1}&=z_{i,t+1}+\delta_{t+1}u_{i,t+1},\\
\lambda_{i,t+1}&=\left[\tilde{\lambda}_{i,t}+\gamma_t(g_{i,t}(x_{i,t})-\beta_t\tilde{\lambda}_{i,t})\right]_+,\label{equ31e}
\end{align}
\end{subequations}
where ${z}_{t}:=col(z_{1,t},\ldots,z_{N,t})$, $\hat{\nabla}_if_{i,t}({z}_{t}):=\frac{n_i}{\delta_t}f_{i,t}(x_{i,t},x_{-i,t})u_{i,t}$, $\hat{\nabla}{g}_{i,t}(z_{i,t}):=\frac{n_i}{\delta_t}g_{i,t}(x_{i,t})u_{i,t}$ $\tilde{\lambda}_{i,t}:=\sum_{j=1}^Na_{ij}\lambda_{j,t}$, $a_{ij}$ is the $(i,j)$th element of adjacency matrix $A$, $\alpha_t,\beta_t,\gamma_t,\eta_t,\delta_t\in[0,1]$ are non-increasing parameters to be determined, and $u_{i,t}\in\mathbb{S}^{n_i}$, $i\in[N]$, are uniformly distributed random vectors.
\end{algorithm}

A distributed bandit online algorithm is given in Algorithm \ref{alg1}, where each player $i\in[N]$ maintains four local variables: the local strategy variable $x_{i,t}\in X_i$, the intermediate local strategy variables $z_{i,t},\tilde{z}_{i,t}\in(1-\eta_t)X_i$, and the local dual variable $\lambda_{i,t}\in\mathbb{R}^m_+$. These variables are updated based on (\ref{equ31b})--(\ref{equ31e}).

The intuition of the update rules in (\ref{equ31b})--(\ref{equ31e}) is explained as follows. A regularized/panelized Lagrangian function at time $t$ associated to each player $i\in[N]$ is defined as
\begin{align}
&{\cal A}_{i,t}(x_{i,t},\lambda_t;x_{-i,t})\nonumber\\
&:=f_{i,t}(x_{i,t},x_{-i,t})+\lambda^{\top}_tg_t(x_t)-\frac{\beta_{t}}{2}\|\lambda_t\|^2,\label{e31}
\end{align}
where $x_{i,t}\in X_i$, and $\lambda_t\in\mathbb{R}^m_+$ is the Lagrange multiplier or dual variable, $\beta_t>0$ is the regularization parameter. It can be known from Lemma 1 in \cite{nedic2009approximate} that there exists a positive constant $\Lambda>0$ to bound the optimal dual variable, that is,
\begin{align}
\|\lambda_t^*\|\leq\Lambda.\label{}
\end{align}

 Based on (\ref{e31}), a mirror-descent-based algorithm is designed as
\begin{align}
x_{i,t+1}
&=\arg\min\limits_{x\in X_i}\{\alpha_{t}\langle x,\nabla_if_{i,t}(x_{i,t},x_{-i,t})\rangle\nonumber\\
&~~~+\alpha_{t}\langle x,(\nabla g_{i,t}(x_{i,t}))^{\top}\lambda_t\rangle+D_{\phi_i}(x,x_{i,t})\},\label{equ33}\\
\lambda_{t+1}&=\left[\lambda_t+\gamma_t(g_t(x_t)-\beta_t\lambda_t)\right]_+,\label{equ34}
\end{align}
where $\alpha_t>0$ and $0<\gamma_t<1$ are stepsizes used in the primal and dual updates, respectively. To update the strategy $x_{i,t}$ of player $i$ needs the information of all other players' strategies, the gradients of $f_{i,t}(x_{i,t},x_{-i,t})$ and $g_{i,t}(x_{i,t})$, the global nonlinear constraint function $g_t(x)$, and the common Lagrange multiplier $\lambda_t$. In contrast, Algorithm \ref{alg1} proposed in this paper is fully distributed and gradient-free.

Time-varying sequences $\{\alpha_{t}\}$, $\{\beta_{t}\}$, $\{\gamma_{t}\}$, $\{\delta_{t}\}$ and $\{\eta_{t}\}$ for all players are assumed to be the same for simplicity in this paper, which in fact can be selected privately for each player. $u_{i,t}$ is randomly chosen following a uniform distribution. Let $\mathfrak{F}_t$ represent the $\sigma$-algebra generated by $(u_{1,t},\ldots,u_{N,t})$, and $\mathcal{F}_t:=\bigcup_{s=1}^t\mathfrak{F}_s$. It can be seen that the random sequences $z_{i,t},x_{i,t-1},\lambda_{i,t},\tilde{\lambda}_{i,t}$, $i\in[N]$ generated by Algorithm \ref{alg1} are independent of $\mathfrak{F}_s$ for $s\geq t$ and depend on $\mathcal{F}_{t-1}$.

In what follows, some necessary lemmas are first presented.
\begin{lemma}\label{lemma3}
Under Assumptions \ref{assump1} and \ref{assump4}, for any $i\in[N]$ and $t\in[T]$, $\lambda_{i,t}$ and $\tilde{\lambda}_{i,t}$ generated by Algorithm \ref{alg1} satisfy
\begin{align}
\|\lambda_{i,t}\|&\leq\frac{B_g}{\beta_t},\label{equ37}\\
\|\tilde{\lambda}_{i,t}\|&\leq\frac{B_g}{\beta_t},\label{equ38}\\
\|\tilde{\lambda}_{i,t}-\overline{\lambda}_t\|&\leq2\sqrt{N}B_g\sum\limits_{s=0}^{t-1}\sigma_m^{s}\gamma_{t-1-s},\label{equ39}\\
\frac{\Lambda_{t+1}}{2\gamma_t}&\leq2NB_g^2\gamma_t+(\overline{\lambda}_t-\lambda)^{\top}g_t(x_t)+\frac{N\beta_t}{2}\|\lambda\|^2\nonumber\\
&~~~+2N\sqrt{N}B_g^2\sum\limits_{s=0}^{t-1}\sigma_m^{s}\gamma_{t-1-s},\label{equ40}
\end{align}
where $\lambda\in\mathbb{R}^m_+$, and
\begin{align*}
\overline{\lambda}_t&:=\frac{1}{N}\sum_{i=1}^N\lambda_{i,t},\\
\Lambda_{t+1}&:=\sum_{i=1}^N\left[\|\lambda_{i,t+1}-\lambda\|^2-(1-\beta_t\gamma_t)\|\lambda_{i,t}-\lambda\|^2\right].
\end{align*}
\end{lemma}

\emph{Proof:} See Appendix \ref{A}. \hfill$\blacksquare$

\begin{lemma}\label{lemma4}
Under Assumptions \ref{assump1}--\ref{assump6}, for all $i\in[N]$, $z_{i,t}$ generated by Algorithm \ref{alg1} satisfies
\begin{align}
&\mu\sum_{t=1}^T\mathbf{E}\left[\|\breve{x}_t^*-z_t\|^2\right]\nonumber\\
&\leq\sum_{t=1}^T\frac{1}{\alpha_t^2}\sum_{i=1}^N\mathbf{E}\left[D_{\phi_i}(\breve{x}_{i,t}^*,z_{i,t})-D_{\phi_i}(\breve{x}_{i,t+1}^*,z_{i,t+1})\right]\nonumber\\
&~~~+\sum_{t=1}^T\frac{K}{\alpha_t^2}\sum_{i=1}^N\|\breve{x}_{i,t+1}^*-\breve{x}_{i,t}^*\|+B_{1,T}+B_{\lambda,T},\label{equ48}
\end{align}
where $\breve{x}_{i,t}^*:=(1-\eta_t)x_{i,t}^*$, $\breve{x}_t^*:=col(\breve{x}_{1,t}^*,\ldots,\breve{x}_{N,t}^*)$, and
\begin{align*}
B_{1,T}&:=2N\sqrt{N}B_xL\sum_{t=1}^T\delta_t+2N(B_xL_g\Lambda+B_g^2)\sum_{t=1}^T\gamma_t\\
&~~+N(2\sqrt{N}B_x^2L+B_xL_f)\sum_{t=1}^T\eta_t+\frac{n^2B_f^2}{\mu_0}\sum_{t=1}^T\frac{\alpha_t}{\delta_t^2}\\
&~~+\frac{n^2B_g^4}{\mu_0}\sum_{t=1}^T\frac{\alpha_t}{\beta_t^2\delta_t^2}+NB_xB_gL_g\sum_{t=1}^T\frac{\eta_t}{\beta_t}\\
&~~+2NB_gL_g\sum_{t=1}^T\frac{\delta_t}{\beta_t}+6N\sqrt{N}B_g^2\sum_{t=1}^T\sum_{s=0}^{t-1}\sigma_m^s\gamma_{t-1-s},
\end{align*}
\begin{align*}
B_{\lambda,T}&:=-\sum_{t=1}^T{\lambda}^{\top}g_{t}(x_{t})+\frac{N}{2}(1+\sum_{t=1}^T\beta_t)\|\lambda\|^2\\
&~~+\sum\limits_{t=1}^T\left(\frac{1}{2\gamma_t}-\frac{1}{2\gamma_{t-1}}-\frac{\beta_t}{2}\right)\sum_{i=1}^N\|\lambda_{i,t}-\lambda\|^2.
\end{align*}
\end{lemma}

\emph{Proof:}  
See Appendix \ref{B}. \hfill$\blacksquare$

\begin{lemma}\label{lemma5}
Under Assumptions \ref{assump1}--\ref{assump6}, for any $i\in[N]$, the expected dynamic regret (\ref{equ3}) and constraint violation (\ref{equ5}) generated by Algorithm \ref{alg1} are bounded by
\begin{align}
\mathbf{E}\left[Reg_i(T)\right]&\leq\frac{L_f}{\sqrt{\mu}}\sqrt{TB_{2,T}+T\Upsilon_{1}+\frac{\sqrt{N}K}{\alpha_T^2}T\Phi_T^*}\nonumber\\
&~~~+B_{3,T},\label{equ67}\\
\mathbf{E}\left\|\left[\sum\limits_{t=1}^Tg_t(x_t)\right]_+\right\|^2&\leq{B}_{2,T}B_{4,T}+B_{4,T}\frac{\sqrt{N}K}{\alpha_T^2}\Phi_T^*\nonumber\\
&~~~+B_{4,T}\Upsilon_{2},
\end{align}
where $\lambda_c:=2B_{4,T}^{-1}\left[\sum_{t=1}^Tg_t(x_t)\right]_+$, and
\begin{align*}
B_{2,T}&:=B_{1,T}+\frac{2NB_xK}{\alpha_{T+1}^2}+NKB_x\sum_{t=1}^T\frac{\eta_t-\eta_{t+1}}{\alpha_t^2},\\
B_{3,T}&:=L_fB_x\sum_{t=1}^T\eta_t+L_f\sum_{t=1}^T\delta_t,\\
B_{4,T}&:=2N(1+\sum_{t=1}^T\beta_t),\\
\Upsilon_{1}&:=\sum\limits_{t=1}^T\left(\frac{1}{2\gamma_t}-\frac{1}{2\gamma_{t-1}}-\frac{\beta_t}{2}\right)\sum_{i=1}^N\|\lambda_{i,t}\|^2,\nonumber\\
\Upsilon_{2}&:=\sum\limits_{t=1}^T\left(\frac{1}{2\gamma_t}-\frac{1}{2\gamma_{t-1}}-\frac{\beta_t}{2}\right)\sum_{i=1}^N\|\lambda_{i,t}-\lambda_c\|^2.
\end{align*}
\end{lemma}

\emph{Proof:} See Appendix \ref{C}. \hfill$\blacksquare$

With the above preparations, it is now ready to present the first main result on Algorithm \ref{alg1} in this paper.
\begin{theorem}\label{thm1}
Under Assumptions \ref{assump1}--\ref{assump6}, let
\begin{align}
\alpha_t=\frac{1}{t^{a_1}}, \beta_t=\frac{1}{t^{a_2}}, \gamma_t=\frac{1}{t^{1-a_2}}, \delta_t=\frac{r_{\min}}{t^{a_3}}, \eta_t=\frac{1}{t^{a_3}},
\end{align}
where the constants $a_1,a_2,a_3\in[0,1]$ satisfy $a_1<0.5$, $a_1-2a_2-2a_3>0$ and $a_2<a_3$, and $r_{\min}:=\min_{i\in[N]}\{r_i\}$.
Then there hold
\begin{align}
\mathbf{E}\left[Reg_{i}(T)\right]&\leq\mathcal{O}\left(T^{\max\{1-\frac{a_1}{2}+a_2+a_3,1+\frac{a_2}{2}-\frac{a_3}{2},\frac{1}{2}+a_1\}}\right)\nonumber\\
&+\mathcal{O}\left(\sqrt{T\log{T}}\right)+\mathcal{O}\left(T^{\frac{1}{2}+a_1}\sqrt{\Phi_T^*}\right),\label{equ75}\\
\mathbf{E}\left[R_g(T)\right]&\leq\mathcal{O}\left(T^{\max\{1-\frac{a_1}{2}+\frac{a_2}{2}+a_3,1-\frac{a_3}{2},\frac{1}{2}+a_1-\frac{a_2}{2}\}}\right)\nonumber\\
&+\mathcal{O}\left(T^{\frac{1}{2}-\frac{a_2}{2}}\sqrt{\log{T}}\right)+\mathcal{O}\left(T^{\frac{1}{2}+a_1-\frac{a_2}{2}}\sqrt{\Phi_T^*}\right).\label{equ76}
\end{align}
\end{theorem}

\emph{Proof:}
For the selected $\beta_t$ and $\gamma_t$, it holds that $\frac{1}{\gamma_t}-\frac{1}{\gamma_{t-1}}-\beta_t\leq0$.
Note that for any $0<a<1$,
\begin{align}
\sum\limits_{t=1}^T\frac{1}{t^a}\leq1+\int_{1}^T\frac{1}{t^{a}}dt\leq1+\frac{T^{1-a}-1}{1-a}\leq\frac{T^{1-a}}{1-a}.
\end{align}
In view of the fact that
\begin{align}
\sum_{t=1}^T\sum\limits_{s=0}^{t-1}\sigma^s_m\gamma_{t-s-1}&=\sum_{t=1}^T\gamma_{t-1}\sum\limits_{s=0}^{T-t}\sigma^s_m\nonumber\\
&\leq\frac{1}{1-\sigma_m}\sum_{t=1}^T\gamma_{t-1},
\end{align}
and similarly,
\begin{align}
\sum_{t=1}^T\sum\limits_{s=0}^{t-1}\sigma^s\alpha_{t-s-1}\leq\frac{1}{1-\sigma}\sum_{t=1}^T\alpha_{t-1}.
\end{align}
In addition, by Lemma 4 in \cite{yi2020distributed_b}, it is obtained that
\begin{align}
\frac{\eta_t-\eta_{t+1}}{\alpha_t^2}&\leq\frac{1}{t},\\
\sum_{t=1}^T\frac{1}{t}&\leq2\log T,~{\rm if}~T>2,
\end{align}
Therefore, $B_{2,T}$, $B_{3,T}$ and $B_{4,T}$ defined in Lemma \ref{lemma5} satisfy
\begin{align}
B_{2,T}&\leq\mathcal{O}(T^{{\max}\{1-a_1+2a_2+2a_3,1+a_2-a_3,2a_1\}})+\mathcal{O}(\log T),\\
B_{3,T}&\leq\mathcal{O}(T^{1-a_3}),\\
B_{4,T}&\leq\mathcal{O}(T^{1-a_2}).
\end{align}
Consequently, invoking Lemma \ref{lemma5} implies (\ref{equ75}) and (\ref{equ76}).\hfill$\blacksquare$

\begin{corollary}\label{corollary1}
Under Assumptions \ref{assump1}--\ref{assump6}, Algorithm \ref{alg1} achieves sublinear ${\bf E}[Reg_i(T)]$ and ${\bf E}[R_g(T)]$ if $\Phi_T^*=\mathcal{O}(T^a)$ for some $a\in[0,1)$.
\end{corollary}

\emph{Proof:}
If $\Phi_T^*=\mathcal{O}(T^a)$, let $a_1\in[0,\frac{1}{2}-\frac{a}{2})$ and choose $a_1,a_2,a_3$ satisfying the conditions in Theorem \ref{thm1}, then $\mathbf{E}[Reg_i(T)]/T\to0,$ and $\mathbf{E}[R_g(T)]/T\to0$ as $T$ goes to infinity. \hfill$\blacksquare$

\begin{corollary}\label{corollary2}
Under Assumptions \ref{assump1}--\ref{assump6}, the same expected dynamic regrets and accumulated constraint violation generated by Algorithm \ref{alg1} are
\begin{align*}
\mathbf{E}\left[Reg_{i}(T)\right]&\leq\mathcal{O}\left(T^{\frac{13}{14}}\right)+\mathcal{O}\left(T^{\frac{13}{14}}\sqrt{\Phi_T^*}\right),\label{}\\
\mathbf{E}\left[R_g(T)\right]&\leq\mathcal{O}\left(T^{\frac{13}{14}}\right)+\mathcal{O}\left(T^{\frac{13}{14}}\sqrt{\Phi_T^*}\right).\label{}
 \end{align*}
\end{corollary}

\emph{Proof}
By Theorem \ref{thm1}, comparing the indexes in (\ref{equ75}) and (\ref{equ76}), the same bounds can be derived by letting $1-\frac{a_1}{2}+a_2+a_3=1+\frac{a_2}{2}-\frac{a_3}{2}=\frac{1}{2}+a_1$, which implies that $a_1=\frac{3}{7}+\frac{4}{7}a_2$ and $a_3=\frac{1}{7}-\frac{1}{7}a_2$. Thus, $1-\frac{a_1}{2}+a_2+a_3=1+\frac{a_2}{2}-\frac{a_3}{2}=\frac{1}{2}+a_1=\frac{13}{14}+\frac{4}{7}a_2$. Note that $a_2\in[0,1]$. The claimed bounds are obtained by making $a_2=0$. \hfill$\blacksquare$


\begin{remark}
To date, a projection-based distributed algorithm was devised in \cite{lu2020online} for online game with time-invariant constraint functions, and then a mirror descent-based distributed algorithm was proposed in \cite{meng2021decentralized} for online game with time-varying constraint functions. However, both works depend on that each player can receive the gradients of its local cost and constraint functions after a strategy profile is determined at each round. In comparison, this paper considers a more challenging scenario, that is, online game with time-varying constraints and one-point bandit feedback, where only function values of cost and constraint functions at the decision vector made by individual agents are revealed gradually. Moreover, the time-varying sequences $\{\alpha_t\}$, $\{\beta_t\}$, $\{\gamma_t\}$, $\{\delta_t\}$ and $\{\eta_t\}$ in this paper are chosen independent of the learning time $T$, which is also a critical issue not addressed in \cite{meng2021decentralized}.
\end{remark}

\section{Distributed Bandit Feedback with Delays}\label{section4}
This section is concerned with the case where the values of $f_{i,t}$ and $g_{i,t}$ are revealed to the player $i$ with time-delays often due to the latency in communication and computation in many real systems, as studied in online optimization \cite{cao2021decentralized}. For such situation, a distributed bandit online algorithm with delays for seeking the variational GNE sequence of the online game is designed as Algorithm \ref{alg2}. The main difference from Algorithm \ref{alg1} is that the delayed information $\hat{\nabla}_if_{i,t-\tau}(z_{t-\tau})$, $\hat{\nabla}{g}_{t-\tau}(z_{i,t-\tau})$ and $g_{i,t-\tau}(x_{i,t-\tau})$ is used at iteration $t+1$ in Algorithm \ref{alg2} rather than the information at time $t$. Here, $\tau>0$ is a constant delay, and is assumed to be the same at agents owing to the common surrounding environment.

\begin{algorithm}[!ht]
\caption{Distributed Bandit Online Primal-Dual Mirror Descent with One-point Delayed Bandit Feedback}\label{alg2}
Each player $i$ maintains vector variables $x_{i,t}\in X_i$, $z_{i,t}\in(1-\eta_t)X_i$, $\tilde{z}_{i,t}\in(1-\eta_t)X_i$  and $\lambda_{i,t}\in\mathbb{R}^{m}_+$ at iteration $t$.

 {\bf Initialization:} For any $i\in[N]$, initialize $z_{i,t},\tilde{z}_{i,t}\in(1-\eta_{t})X_i$ arbitrarily, $x_{i,t}=z_{i,t}+\delta_tu_{i,t}$, and $\lambda_{i,t}={\bf0}_m$, $t\in[1,\tau+1]$.

{\bf Iteration:}
At each $t>\tau$, every player $i$ receives the delayed function values $f_{i,t-\tau}(x_{t-\tau})$ and $g_{i,t-\tau}(x_{i,t-\tau})$, and performs the following update:
\begin{subequations}
\begin{align}
\tilde{z}_{i,t+1}&=\arg\min\limits_{z\in(1-\eta_t)X_i}\{\alpha_t\langle{z},\hat{\nabla}_if_{i,t-\tau}(z_{t-\tau})\rangle\nonumber\\
&~~~~~~~~~~~~~~+\alpha_t\langle z,(\hat{\nabla} g_{i,t-\tau}(z_{i,t-\tau}))^{\top}\tilde{\lambda}_{i,t}\rangle\nonumber\\
&~~~~~~~~~~~~~~+D_{\phi_i}(z,z_{i,t-\tau})\},\label{e54b}\\
z_{i,t+1}&=(1-\alpha_t)z_{i,t-\tau}+\alpha_t\tilde{z}_{i,t+1},\label{}\\
x_{i,t+1}&=z_{i,t+1}+\delta_{t+1}u_{i,t+1},\\
\lambda_{i,t+1}&=\left[\tilde{\lambda}_{i,t}+\gamma_t(g_{i,t-\tau}(x_{i,t-\tau})-\beta_t\tilde{\lambda}_{i,t})\right]_+,\label{}
\end{align}
\end{subequations}
where $z_{t}:=col(z_{1,t},\ldots,z_{N,t})$, $\tilde{\lambda}_{i,t}:=\sum_{j=1}^Na_{ij}\lambda_{j,t}$, $a_{ij}$ is the $(i,j)$th element of adjacency matrix $A$, $\alpha_t,\beta_t,\gamma_t,\eta_t,\delta_t\in[0,1]$ are non-increasing sequences to be determined, and $u_{i,t}\in\mathbb{S}^{n_i}$, $i\in[N]$, are uniformly distributed random vectors.
\end{algorithm}

Similar to the results on Algorithm \ref{alg1}, the following results on Algorithm \ref{alg2} are obtained.
\begin{lemma}\label{}
Under Assumptions \ref{assump1} and \ref{assump4}, for any $i\in[N]$ and $t>\tau$, $\lambda_{i,t}$ and $\tilde{\lambda}_{i,t}$ generated by Algorithm \ref{alg2} satisfy
\begin{align}
\|\lambda_{i,t}\|&\leq\frac{B_g}{\beta_t},\label{}\\
\|\tilde{\lambda}_{i,t}\|&\leq\frac{B_g}{\beta_t},\label{}\\
\|\tilde{\lambda}_{i,t}-\overline{\lambda}_t\|&\leq2\sqrt{N}B_g\sum\limits_{s=0}^{t-1}\sigma_m^{s}\gamma_{t-1-s},\label{}\\
\frac{\Lambda_{t+1}}{2\gamma_t}&\leq2NB_g^2\gamma_t+(\overline{\lambda}_t-\lambda)^{\top}g_{t-\tau}(x_{t-\tau})+\frac{N\beta_t}{2}\|\lambda\|^2\nonumber\\
&~~~+2N\sqrt{N}B_g^2\sum\limits_{s=0}^{t-1}\sigma_m^{s}\gamma_{t-1-s},\label{}
\end{align} where $\lambda\in\mathbb{R}^m_+$, $\overline{\lambda}_t:=\frac{1}{N}\sum_{i=1}^N\lambda_{i,t}$, and
\begin{align*}
\Lambda_{t+1}&:=\sum_{i=1}^N\left[\|\lambda_{i,t+1}-\lambda\|^2-(1-\beta_t\gamma_t)\|\lambda_{i,t}-\lambda\|^2\right].
\end{align*}
\end{lemma}

\emph{Proof:}
The proof is the same as that of Lemma \ref{lemma3}, which is thus omitted here. \hfill$\blacksquare$

\begin{lemma}\label{lemma8}
Under Assumptions \ref{assump1}--\ref{assump6}, for all $i\in[N]$, $z_{i,t}$ generated by Algorithm \ref{alg2} satisfies
\begin{align}
&\mu\sum_{t=1}^{T}\mathbf{E}\left[\|\breve{x}_{t}^*-z_{t}\|^2\right]\nonumber\\
&=\mu\sum_{t=\tau+1}^{T+\tau}\mathbf{E}\left[\|\breve{x}_{t-\tau}^*-z_{t-\tau}\|^2\right]\nonumber\\
&\leq\sum_{t=\tau+1}^{T+\tau}\frac{1}{\alpha_t^2}\sum_{i=1}^N\mathbf{E}\left[D_{\phi_i}(\breve{x}_{i,t-\tau}^*,z_{i,t-\tau})-D_{\phi_i}(\breve{x}_{i,t+1}^*,z_{i,t+1})\right]\nonumber\\
&~~~+\sum_{t=\tau+1}^{T+\tau}\frac{K}{\alpha_t^2}\sum_{i=1}^N\|\breve{x}_{i,t+1}^*-\breve{x}_{i,t-\tau}^*\|+C_{1,T}+C_{\lambda,T},\label{equ91}
\end{align}
where
\begin{align*}
C_{1,T}&:=2N\sqrt{N}B_xL\sum_{t=\tau+1}^{T+\tau}\delta_{t-\tau}+2NB_gL_g\sum_{t=\tau+1}^{T+\tau}\frac{\delta_{t-\tau}}{\beta_t}\\
&~~~~+2N(B_xL_g\Lambda+B_g^2)\sum_{t=\tau+1}^{T+\tau}\gamma_t+\frac{n^2B_f^2}{\mu_0}\sum_{t=\tau+1}^{T+\tau}\frac{\alpha_t}{\delta_{t-\tau}^2}\\
&~~~~+N(2\sqrt{N}B_x^2L+B_xL_f)\sum_{t=\tau+1}^{T+\tau}\eta_{t-\tau}\\
&~~~~+\frac{n^2B_g^4}{\mu_0}\sum_{t=\tau+1}^{T+\tau}\frac{\alpha_t}{\beta_t^2\delta_{t-\tau}^2}+NB_xB_gL_g\sum_{t=\tau+1}^{T+\tau}\frac{\eta_{t-\tau}}{\beta_{t}}\\
&~~~~+6N\sqrt{N}B_g^2\sum_{t=\tau+1}^{T+\tau}\sum_{s=0}^{t-1}\sigma_m^s\gamma_{t-1-s},\\
C_{\lambda,T}&:=-\sum_{t=\tau+1}^{T+\tau}{\lambda}^{\top}g_{t-\tau}(x_{t-\tau})+\frac{N}{2}(1+\sum_{t=\tau+1}^{T+\tau}\beta_t)\|\lambda\|^2\\
&~~~~+\sum\limits_{t=\tau+1}^{T+\tau}\left(\frac{1}{2\gamma_t}-\frac{1}{2\gamma_{t-1}}-\frac{\beta_t}{2}\right)\sum_{i=1}^N\|\lambda_{i,t}-\lambda\|^2,
\end{align*}
\end{lemma}

\emph{Proof:}
See Appendix \ref{E}. \hfill$\blacksquare$

\begin{lemma}\label{lemma9}
Under Assumptions \ref{assump1}--\ref{assump6}, for any $i\in[N]$, the expected dynamic regret (\ref{equ3}) and constraint violation (\ref{equ5}) generated by Algorithm \ref{alg2} is bounded by
\begin{align}
&\mathbf{E}\left[Reg_i(T)\right]\nonumber\\
&\leq\frac{L_f}{\sqrt{\mu}}\sqrt{TC_{2,T}+T\Upsilon_{\tau,1}+\frac{\sqrt{N}K}{\alpha_{T+\tau}^2}(\tau+1){T}\Phi_{T+\tau}^*}+C_{3,T},\label{equ92}\\
&\mathbf{E}\left\|\left[\sum\limits_{t=1}^Tg_t(x_t)\right]_+\right\|^2\nonumber\\
&\leq{C}_{2,T}C_{4,T}+C_{4,T}\frac{\sqrt{N}K}{\alpha_{T+\tau}^2} (\tau+1)\Phi_{T+\tau}^*+C_{4,T}\Upsilon_{\tau,2},
\end{align}
where $\lambda_{\tau}:=2C_{4,T}^{-1}\left[\sum_{t=1}^Tg_t(x_t)\right]_+$,  and
\begin{align*}
C_{2,T}&:=C_{1,T}+\frac{2NB_xK(\tau+1)}{\alpha_{T+\tau}^2}\\
&~~~~+NKB_x(\tau+1)\sum_{t=\tau+1}^{T+\tau}\frac{\eta_{t-\tau}-\eta_{t+1-\tau}}{\alpha_t^2},\\
C_{3,T}&:=L_fB_x\sum_{t=1}^T\eta_t+L_f\sum_{t=1}^T\delta_t,\\
C_{4,T}&:=2N(1+\sum_{t=\tau+1}^{T+\tau}\beta_t),\\
\Upsilon_{\tau,1}&:=\sum\limits_{t=\tau+1}^{T+\tau}\left(\frac{1}{2\gamma_t}-\frac{1}{2\gamma_{t-1}}-\frac{\beta_t}{2}\right)\sum_{i=1}^N\|\lambda_{i,t}\|^2,\\
\Upsilon_{\tau,2}&:=\sum\limits_{t=\tau+1}^{T+\tau}\left(\frac{1}{2\gamma_t}-\frac{1}{2\gamma_{t-1}}-\frac{\beta_t}{2}\right)\sum_{i=1}^N\|\lambda_{i,t}-\lambda_{\tau}\|^2.
\end{align*}
\end{lemma}

\emph{Proof:}
See Appendix \ref{F}. \hfill$\blacksquare$

Equipped with the above preparations, we are now ready to present the second main result of this paper.
\begin{theorem}\label{thm2}
Under Assumptions \ref{assump1}--\ref{assump6}, choose $\alpha_t=\beta_t=\gamma_t=1$ for $t\in[0,\tau]$ and
\begin{align}
&\alpha_t=\frac{1}{(t-\tau)^{a_1}}, \beta_t=\frac{1}{(t-\tau)^{a_2}}, \gamma_t=\frac{1}{(t-\tau)^{1-a_2}},~t>\tau\nonumber\\ &\delta_t=\frac{r_{\min}}{t^{a_3}}, \eta_t=\frac{1}{t^{a_3}}, ~t\geq1,
\end{align}
where the constants $a_1,a_2,a_3\in[0,1]$ satisfy $a_1-2a_2-2a_3>0$ and $a_2<a_3$, and $r_{\min}:=\min_{i\in[N]}\{r_i\}$.
Then for each $i\in[N]$, there hold
\begin{align}
\mathbf{E}\left[Reg_{i}(T)\right]&\leq\mathcal{O}\left(T^{\max\{1-\frac{a_1}{2}+a_2+a_3,1+\frac{a_2}{2}-\frac{a_3}{2}\}}\right)\nonumber\\
&~~~+\mathcal{O}(\sqrt{\tau+1}{T}^{\frac{1}{2}+a_1})+\mathcal{O}\left(\sqrt{(\tau+1){T}\log{T}}\right)\nonumber\\
&~~~+\mathcal{O}\left(T^{\frac{1}{2}+a_1}\sqrt{(\tau+1)\Phi_{T+\tau}^*}\right),\label{}\\
\mathbf{E}\left[R_g(T)\right]&\leq\mathcal{O}\left(T^{\max\{1-\frac{a_1}{2}+\frac{a_2}{2}+a_3,1-\frac{a_3}{2}\}}\right)\nonumber\\
&~~~+\mathcal{O}\left(\sqrt{\tau+1}{T}^{\frac{1}{2}+a_1-\frac{a_2}{2}}\right)\nonumber\\
&~~~+\mathcal{O}\left({T}^{\frac{1}{2}-\frac{a_2}{2}}\sqrt{(\tau+1)\log{T}}\right)\nonumber\\
&~~~+\mathcal{O}\left({T}^{\frac{1}{2}+a_1-\frac{a_2}{2}}\sqrt{(\tau+1){\Phi_{T+\tau}^*}}\right).\label{}
\end{align}
\end{theorem}

\emph{Proof:}
Similar to the proof of Theorem \ref{thm1}, it can be derived that
\begin{align}
C_{1,T}&\leq\mathcal{O}(T^{\max\{T^{1-a_1+2a_2+2a_3,1+a_2-a_3}\}}),\\
C_{2,T}&\leq\mathcal{O}(T^{\max\{T^{1-a_1+2a_2+2a_3,1+a_2-a_3}\}})\nonumber\\
&~~~+\mathcal{O}((\tau+1){T}^{2a_1}+(\tau+1)\log{T}),\\
C_{3,T}&\leq\mathcal{O}(T^{1-a_3}),\\
C_{4,T}&\leq\mathcal{O}(T^{1-a_2}).
\end{align}
Then, the proof is completed via a simple computation. \hfill$\blacksquare$

Similar to Corollaries \ref{corollary1} and \ref{corollary2}, the following results are obtained.
\begin{corollary}
Under Assumptions \ref{assump1}--\ref{assump6}, there hold for Algorithm \ref{alg2} that the expected dynamic regrets and accumulated constraint violation are sublinearly bounded if both $(\tau+1)\log{T}$ and $(\tau+1){\Phi_{T+\tau}^*}$ are sublinear with respect to $T$.
\end{corollary}
\begin{corollary}
Under Assumptions \ref{assump1}--\ref{assump6}, the same bounds of the expected dynamic regrets and accumulated constraint violation generated by Algorithm \ref{alg2} are
\begin{align*}
\mathbf{E}\left[Reg_{i}(T)\right]&\leq\mathcal{O}\left(\sqrt{\tau+1}T^{\frac{13}{14}}\right)+\mathcal{O}\left(T^{\frac{13}{14}}\sqrt{(\tau+1)\Phi_{T+\tau}^*}\right),\label{}\\
\mathbf{E}\left[R_g(T)\right]&\leq\mathcal{O}\left(\sqrt{\tau+1}T^{\frac{13}{14}}\right)+\mathcal{O}\left(T^{\frac{13}{14}}\sqrt{(\tau+1)\Phi_{T+\tau}^*}\right).\label{}
 \end{align*}
\end{corollary}
\begin{remark}
The bounds of the expected regrets and constraint violation are related to the time-delay $\tau$ and $\Phi_{T+\tau}^*$, indicating that to learn the GNE within $T$ requires that the path variation of GNE from time $t=1$ to time $t=T+\tau$ increases sublinearly. In fact, if we consider static regret, where $x_t^*=x^*$ for all $t$, then the term $\Phi_{T+\tau}^*$ is 0 and thus will disappear. Note that all the agents are assumed to undergo the same delay when computing their function values in this paper without considering communication delays when transmitting information among agents. In this respect, it is also interesting to see if the same results can be affirmed in the case with different delays for individual agents and/or communication delays, which is left as one of future works.
\end{remark}

\section{Numerical Examples}\label{section5}
\begin{figure}[!ht]
 \centering
  \includegraphics[width=3in]{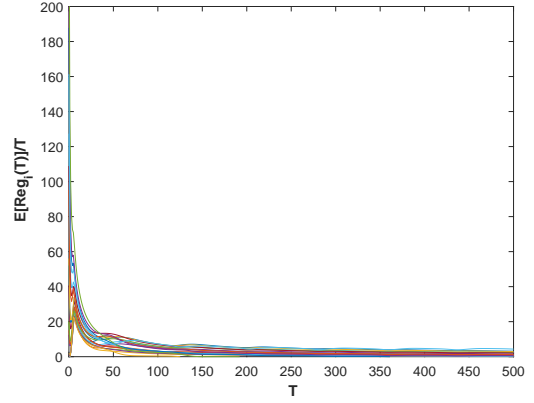}
  \caption{The expected average regrets $\mathbf{E}[Reg_i(T)/T]$, $i\in[N]$ in Algorithm \ref{alg1}.}\label{fig1}
\end{figure}

\begin{figure}[!ht]
 \centering
  \includegraphics[width=3in]{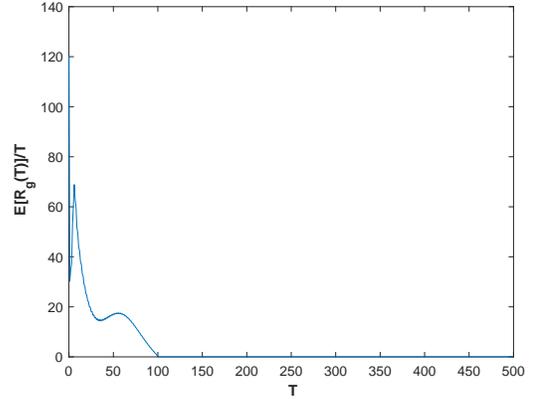}
  \caption{The expected average violation $\mathbf{E}[R_g(T)/T]$, $i\in[N]$ in Algorithm \ref{alg1}.}\label{fig2}
\end{figure}

\begin{figure}[!ht]
 \centering
  \includegraphics[width=3in]{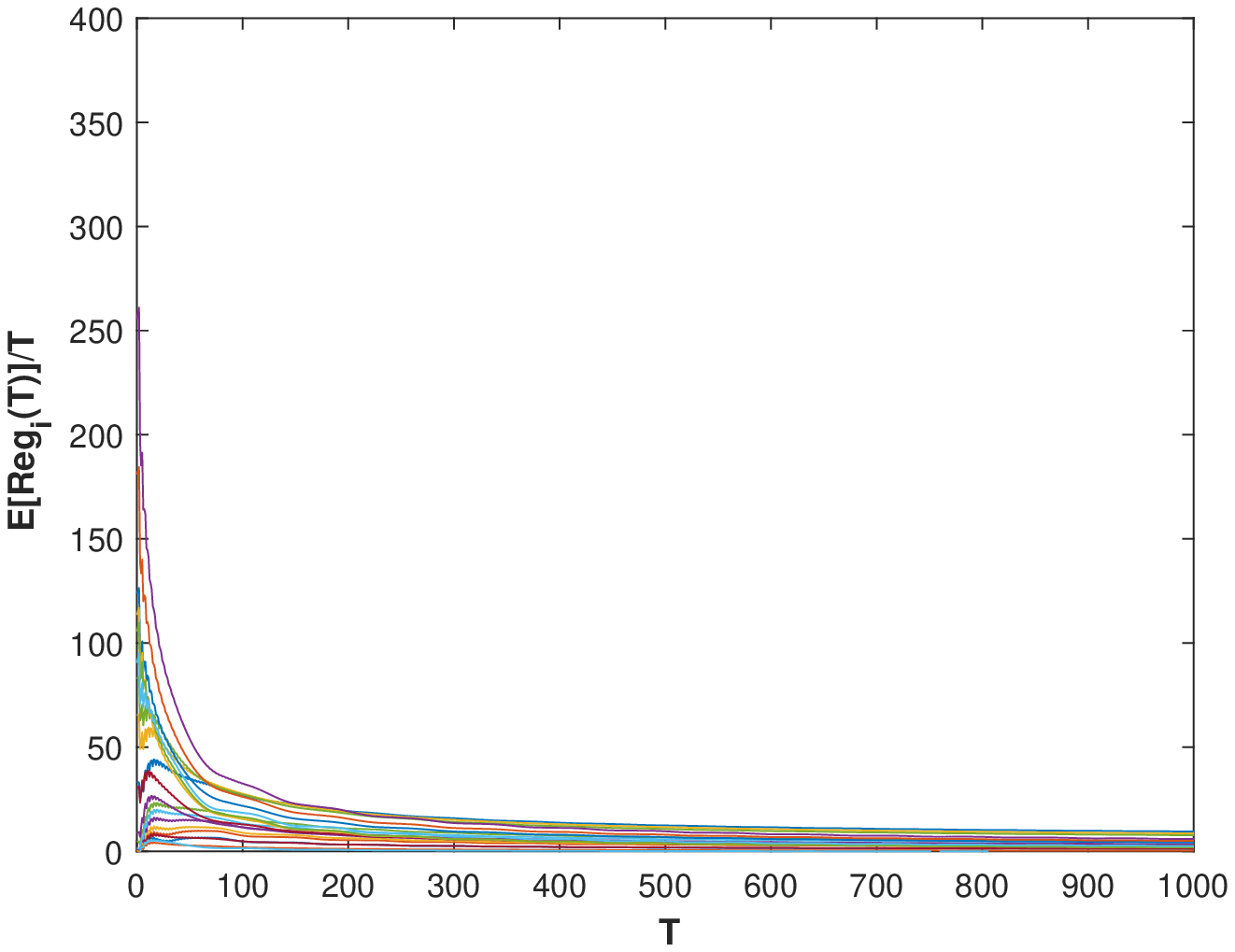}
  \caption{The expected regrets ${\bf E}[Reg_i(T)]/T$, $i\in[N]$ in Algorithm \ref{alg2}.}\label{fig3}
\end{figure}

\begin{figure}[!ht]
 \centering
  \includegraphics[width=3in]{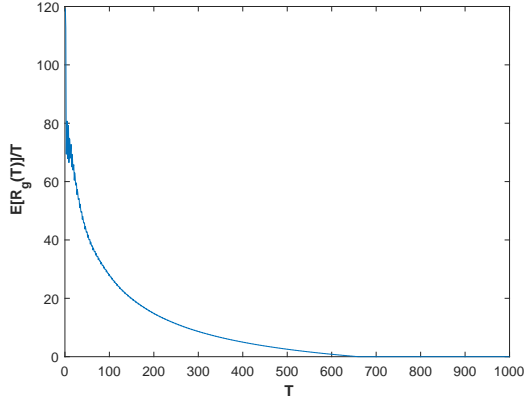}
  \caption{The expected regrets ${\bf E}[R_g(T)]/T$, $i\in[N]$ in Algorithm \ref{alg2}.}\label{fig4}
\end{figure}

An online Nash-Cournot game $\Gamma([N],X_t,f_t)$ is adopted to illustrate the feasibility of the proposed algorithms. Assume that there are $N=20$ firms in the Nash-Cournot game, competing on the amount of product that they will produce subject to time-varying production constraints and market capacity constraints. The quality produced by firm $i$ at time $t$ is denoted as $x_{i,t}\in{X}_i=[0,30]$ and the cost function of firm $i$ is written as $f_{i,t}(x_{i,t},x_{-i,t})=p_{i,t}(x_{i,t})-d_{i,t}(x_{t})$ with the production cost being $p_{i,t}(x_{i,t})=x_{i,t}(\sin(t/12)+1)$ and the demand price of firm $i$ being $d_{i,t}(x_{t})=21+i/9-0.5i\sin(t/12)-\sum_{j=1}^Nx_{j,t}$. The market capacity constraint is considered as the shared inequality constraint $\sum_{i=1}^Nx_{i,t}\leq\sum_{i=1}^Nl_{i,t}$, where $l_{i,t}=10+\sin(t/12)$ is the local bound only available to firm $i$. If considering the offline and centralized setting, the time-varying GNE can be computed as $x_{i,t}^*=P_{X_i}(\xi_{i,t})$, where $\xi_{i,t}:=\frac{1}{9}(i-10)+\frac{1}{2}(10-i)\sin\frac{t}{12}$. In the distributed and bandit feedback setting, choose $a_1=0.45,a_2=0.1,a_3=0.11$, and set initial states $z_{i,1}\in(1-\eta_1){X}_{i}$ randomly, and $\lambda_{i,1}={\bf0}_m$, then by Algorithm \ref{alg1}, the averages of dynamic regrets and constraint violation are shown in Figs. \ref{fig1} and \ref{fig2}, respectively, from which one can see that Algorithm \ref{alg1} achieves sublinear expected regrets and constraint violation. Furthermore, if feedback delays are considered and let $\tau=1$, then the averages of dynamic regrets and constraint violation are shown in Figs. \ref{fig3} and \ref{fig4}, respectively, from which it can be seen that Algorithm \ref{alg2} achieves sublinear expected regrets and constraint violation. Moreover, feedback delay  makes the convergence of the averages of dynamic regrets and constraint violation generated by Algorithm \ref{alg2} a little bit slower than that generated by Algorithm \ref{alg1}.

\section{Conclusion}\label{section6}
In this paper, distributed online learning of GNEs in online game was studied based on bandit feedback and mirror descent. Each player in online game aims to minimize its own self-interested time-varying cost function subject to a time-varying coupled inequality constraints with no information of gradients of cost and constraint functions. To this end, distributed online algorithms based on mirror descent and one-point bandit feedback were designed for online game without and with feedback delays, respectively, which can achieve sublinear dynamic regrets and accumulated constraint violation if some quantities are sublinear. Future research directions can be placed on developing a method to improve the regret bounds, to consider the scenario of bandit feedback with different time-delays for different agents and to study the case with communication delays when transmitting information among agents besides delayed bandit feedback.
\begin{appendix}
\subsection{Proof of Lemma \ref{lemma3}}\label{A}
(\ref{equ37}) and (\ref{equ38}) can be easily proved by mathematical induction.

Let us now prove (\ref{equ39}). From (\ref{equ31e}), it can be obtained that
\begin{align}\label{equ41}
\lambda_{i,t+1}=\tilde{\lambda}_{i,t}+\varepsilon_{i,t},
\end{align}
where
$\varepsilon_{i,t}
:=\left[\tilde{\lambda}_{i,t}+\gamma_t(g_{i,t}(x_{i,t})-\beta_t\tilde{\lambda}_{i,t})\right]_+-\tilde{\lambda}_{i,t}$, $i\in[N]$.
Let $\hat{\lambda}_{i,t}:=\lambda_{i,t}-\overline{\lambda}_{t}$, $i\in[N]$, $\hat{\lambda}_t:=col(\hat{\lambda}_{1,t},\ldots,\hat{\lambda}_{N,t})$, and $\varepsilon_t=col(\varepsilon_{1,t},\ldots,\varepsilon_{N,t})$.  Then, the iteration of $\hat{\lambda}_t$ is given as
\begin{align}
\hat{\lambda}_{t+1}&=\left((A-\frac{1}{N}{\bf1}_N{\bf1}_N^{\top})\otimes I_m\right)\hat{\lambda}_t\nonumber\\
&~~~+\left((I-\frac{1}{N}{\bf1}_N{\bf1}_N^{\top})\otimes I_m\right)\varepsilon_t. \label{equ42}
\end{align}
Under Assumption \ref{assump4}, one has $0\leq\sigma_m<1$, where $\sigma_m=\|A-\frac{1}{N}{\bf1}_N{\bf1}_N^{\top}\|$. Taking norm on both sides of (\ref{equ42}) leads to that
\begin{align}
\|\hat{\lambda}_{t+1}\|&\leq\sigma_m\|\hat{\lambda}_{t}\|+\|\varepsilon_t\|.
\end{align}
Note that
\begin{align}
\|\varepsilon_{i,t}\|
&=\left\|\left[\tilde{\lambda}_{i,t}+\gamma_t(g_{i,t}(x_{i,t})-\beta_t\tilde{\lambda}_{i,t})\right]_+-\tilde{\lambda}_{i,t}\right\|\nonumber\\
&\leq\gamma_t\left\|g_{i,t}(x_{i,t})-\beta_t\tilde{\lambda}_{i,t}\right\|\nonumber\\
&\leq \gamma_t\|g_{i,t}(x_{i,t})\|+\gamma_t\beta_t\|\tilde{\lambda}_{i,t}\|\nonumber\\
&\leq 2B_g\gamma_t,\label{}
\end{align}
where the first inequality relies on the inequality $\|[x]_+-[y]_+\|\leq\|x-y\|$ for any $x,y\in\mathbb{R}^m$, and the third inequality is derived following (\ref{equ7}) and (\ref{equ38}).
It yields that
\begin{align}
\|\hat{\lambda}_{t+1}\|\leq\sigma_m\|\hat{\lambda}_{t}\|+2\sqrt{N}B_g\gamma_t.
\end{align}
Note that $\hat{\lambda}_{1}={\bf0}_{Nm}$. Thus, (\ref{equ39}) is proved.

As for (\ref{equ40}), it can be derived from (\ref{equ31e}) that for $\lambda\in\mathbb{R}^m_+$,
\begin{align}
\|\lambda_{i,t+1}-\lambda\|^2
&=\left\|\left[\tilde{\lambda}_{i,t}+\gamma_t(g_{i,t}(x_{i,t})-\beta_t\tilde{\lambda}_{i,t})\right]_+-\lambda\right\|^2\nonumber\\
&\leq\left\|\tilde{\lambda}_{i,t}+\gamma_t(g_{i,t}(x_{i,t})-\beta_t\tilde{\lambda}_{i,t})-\lambda\right\|^2\nonumber\\
&=\|\tilde{\lambda}_{i,t}-\lambda\|^2+\gamma_t^2\|g_{i,t}(x_{i,t})-\beta_t\tilde{\lambda}_{i,t}\|^2\nonumber\\
&~~~+2\gamma_t(\tilde{\lambda}_{i,t}-\lambda)^{\top}(g_{i,t}(x_{i,t})-\beta_t\tilde{\lambda}_{i,t})\nonumber\\
&\leq\sum_{i=1}^Na_{ij}\|\lambda_{j,t}-\lambda\|^2+4B_g^2\gamma_t^2\nonumber\\
&~~~+2\gamma_t(\tilde{\lambda}_{i,t}-\lambda)^{\top}(g_{i,t}(x_{i,t})-\beta_t\tilde{\lambda}_{i,t}),\label{equ46}
\end{align}
where the last inequality is based on (\ref{equ7}) and (\ref{equ38}).

For the last term in (\ref{equ46}),
\begin{align}
&2\gamma_t(\tilde{\lambda}_{i,t}-\lambda)^{\top}(g_{i,t}(x_{i,t})-\beta_t\tilde{\lambda}_{i,t})\nonumber\\
&=2\gamma_t(\tilde{\lambda}_{i,t}-\overline{\lambda}_t)^{\top}g_{i,t}(x_{i,t})+2\gamma_t(\overline{\lambda}_t-\lambda)^{\top}g_{i,t}(x_{i,t})\nonumber\\
&~~~-2\beta_t\gamma_t(\tilde{\lambda}_{i,t}-\lambda)^{\top}\tilde{\lambda}_{i,t}\nonumber\\
&\leq2\gamma_t\|\tilde{\lambda}_{i,t}-\overline{\lambda}_t\|\|g_{i,t}(x_{i,t})\|+2\gamma_t(\overline{\lambda}_t-\lambda)^{\top}g_{i,t}(x_{i,t})\nonumber\\
&~~~+\beta_t\gamma_t(\|\lambda\|^2-\|\tilde{\lambda}_{i,t}-\lambda\|^2)\nonumber\\
&\leq4\sqrt{N}B_g^2\gamma_t\sum\limits_{s=0}^{t-1}\sigma_m^{s}\gamma_{t-1-s}+2\gamma_t(\overline{\lambda}_t-\lambda)^{\top}g_{i,t}(x_{i,t})\nonumber\\
&~~~+\beta_t\gamma_t(\|\lambda\|^2-\|\tilde{\lambda}_{i,t}-\lambda\|^2). \label{equ47}
\end{align}
Substituting (\ref{equ47}) into (\ref{equ46}), summing over $i\in[N]$ on both sides of (\ref{equ46}), and then a simple computation leads to (\ref{equ40}). \hfill$\blacksquare$

\subsection{Proof of Lemma \ref{lemma4}}\label{B}
For any $i\in[N]$, it can be derived based on the optimality of $\tilde{z}_{i,t+1}$ in (\ref{equ31b}) that for any $z\in(1-\eta_t)X_i$,
\begin{align}
&\left\langle\tilde{z}_{i,t+1}-z,\alpha_t(\hat{\nabla}_if_{i,t}({z}_{t})+(\hat{\nabla} g_{i,t}(z_{i,t}))^{\top}\tilde{\lambda}_{i,t})\right\rangle\nonumber\\
&+\left\langle\tilde{z}_{i,t+1}-z,\nabla\phi_i(\tilde{z}_{i,t+1})-\nabla\phi_i({z}_{i,t})\right\rangle\leq 0,\label{equ50}
\end{align}
where $\frac{\partial D_{\phi_i}(z,z_{i,t})}{\partial z}=\nabla\phi_i({z})-\nabla\phi_i({z}_{i,t})$ for $z\in X_i$ is adopted. Taking $z=\breve{x}^*_{i,t}$, from (\ref{equ50}), one has that
\begin{align}
&\alpha_t\left\langle\tilde{z}_{i,t+1}-\breve{x}_{i,t}^*,\hat{\nabla}_if_{i,t}({z}_{t})+(\hat{\nabla}g_{i,t}(z_{i,t}))^{\top}\tilde{\lambda}_{i,t}\right\rangle\nonumber\\
&\leq\left\langle\breve{x}_{i,t}^*-\tilde{z}_{i,t+1},\nabla\phi_i(\tilde{z}_{i,t+1})-\nabla\phi_i({z}_{i,t})\right\rangle\nonumber\\
&=D_{\phi_i}(\breve{x}_{i,t}^*,z_{i,t})-D_{\phi_i}(\breve{x}_{i,t}^*,\tilde{z}_{i,t+1})-D_{\phi_i}(\tilde{z}_{i,t+1},z_{i,t})\nonumber\\
&\leq D_{\phi_i}(\breve{x}_{i,t}^*,z_{i,t})-D_{\phi_i}(\breve{x}_{i,t}^*,\tilde{z}_{i,t+1})-\frac{\mu_0}{2}\|\tilde{z}_{i,t+1}-z_{i,t}\|^2,\label{equ51}
\end{align}
where the equality is based on (\ref{equ18}) and the second inequality is derived via (\ref{equ17}). Therefore,
\begin{align}
&D_{\phi_i}(\breve{x}_{i,t}^*,\tilde{z}_{i,t+1})\nonumber\\
&\leq\alpha_t\left\langle\breve{x}_{i,t}^*-\tilde{z}_{i,t+1},(\hat{\nabla}_if_{i,t}({z}_{t})\right\rangle\nonumber\\
&~~~+\alpha_t\left\langle\breve{x}_{i,t}^*-\tilde{z}_{i,t+1},(\hat{\nabla}g_{i,t}(z_{i,t}))^{\top}\tilde{\lambda}_{i,t})\right\rangle\nonumber\\
&~~~+ D_{\phi_i}(\breve{x}_{i,t}^*,z_{i,t})-\frac{\mu_0}{2}\|\tilde{z}_{i,t+1}-z_{i,t}\|^2.\label{equ52}
\end{align}

Taking expectation with respect to $u_{i,t}$ on both sides of (\ref{equ52}) yields that
\begin{align}
&\mathbf{E}_{u_{i,t}}\left[D_{\phi_i}(\breve{x}_{i,t}^*,\tilde{z}_{i,t+1})\right]\nonumber\\
&\leq\alpha_t\mathbf{E}_{u_{i,t}}\left[\left\langle\breve{x}_{i,t}^*-\tilde{z}_{i,t+1},\hat{\nabla}_if_{i,t}({z}_{t})\right\rangle\right]\nonumber\\
&~~~+\alpha_t\mathbf{E}_{u_{i,t}}\left[\left\langle\breve{x}_{i,t}^*-\tilde{z}_{i,t+1},(\hat{\nabla}g_{i,t}(z_{i,t}))^{\top}\tilde{\lambda}_{i,t}\right\rangle\right]\nonumber\\
&~~~+\mathbf{E}_{u_{i,t}}\left[D_{\phi_i}(\breve{x}_{i,t}^*,z_{i,t})\right]-\frac{\mu_0}{2}\mathbf{E}_{u_{i,t}}\left[\|\tilde{z}_{i,t+1}-z_{i,t}\|^2\right]\nonumber\\
&=\alpha_t\left\langle\breve{x}_{i,t}^*-{z}_{i,t},{\nabla}_i\hat{f}_{i,t}({z}_{t})\right\rangle\nonumber\\
&~~~+\alpha_t\mathbf{E}_{u_{i,t}}\left[\left\langle{z}_{i,t}-\tilde{z}_{i,t+1},\hat{\nabla}_if_{i,t}({z}_{t})\right\rangle\right]\nonumber\\
&~~~+\alpha_t\mathbf{E}_{u_{i,t}}\left[\left\langle\breve{x}_{i,t}^*-\tilde{z}_{i,t+1},(\hat{\nabla}g_{i,t}(z_{i,t}))^{\top}\tilde{\lambda}_{i,t}\right\rangle\right]\nonumber\\
&~~~+D_{\phi_i}(\breve{x}_{i,t}^*,z_{i,t})-\frac{\mu_0}{2}\mathbf{E}_{u_{i,t}}\left[\|\tilde{z}_{i,t+1}-z_{i,t}\|^2\right],\label{equ53}
\end{align}
where the equality holds since $z_{i,t}$ is independent of $u_{i,t}$ and $\mathbf{E}_{u_{i,t}}[\hat{\nabla}_if_{i,t}(z_t)]=\nabla_i\hat{f}_{i,t}(z_t)$ by Lemma \ref{lemma1} with $\hat{f}_{i,t}(z_t):=\mathbb{E}_{v_{i,t}\in\mathbb{B}^{n_i}}[f_{i,t}(z_{i,t}+\delta_tv_{i,t},x_{-i,t})]$.
Next, we discuss the bound of each term on the right-hand side of (\ref{equ53}).

For the first term on the right-hand side of (\ref{equ53}), one has that
\begin{align}
&\alpha_t\left\langle\breve{x}_{i,t}^*-{z}_{i,t},\nabla_i\hat{f}_{i,t}(z_t)\right\rangle\nonumber\\
&=\alpha_t\left\langle\breve{x}_{i,t}^*-{z}_{i,t},\nabla_if_{i,t}(z_t)-\nabla_if_{i,t}(\breve{x}_{t}^*)\right\rangle\nonumber\\
&~~~+\alpha_t\left\langle\breve{x}_{i,t}^*-{z}_{i,t},\nabla_i\hat{f}_{i,t}(z_t)-\nabla_if_{i,t}(z_t)\right\rangle\nonumber\\
&~~~+\alpha_t\left\langle\breve{x}_{i,t}^*-{z}_{i,t},\nabla_if_{i,t}(\breve{x}_{t}^*)-{\nabla}_if_{i,t}(x_t^*)\right\rangle\nonumber\\
&~~~+\alpha_t\left\langle\breve{x}_{i,t}^*-{z}_{i,t},{\nabla}_if_{i,t}(x_t^*)\right\rangle.\label{equ54}
\end{align}
For the second term on the right-hand side of (\ref{equ54}), it is obtained that
\begin{align}
&\alpha_t\left\langle\breve{x}_{i,t}^*-{z}_{i,t},\nabla_i\hat{f}_{i,t}(z_t)-\nabla_if_{i,t}(z_t)\right\rangle\nonumber\\
&=\alpha_t\mathbf{E}_{v_{i,t}\in\mathbb{B}^{n_i}}[\langle\breve{x}_{i,t}^*-{z}_{i,t},\nabla_i{f}_{i,t}(z_{i,t}+\delta_tv_{i,t},x_{-i,t})\nonumber\\
&~~~~~-\nabla_if_{i,t}(z_t)\rangle]\nonumber\\
&\leq\alpha_t\mathbf{E}_{v_{i,t}\in\mathbb{B}^{n_i}}[2B_x\|\nabla_i{f}_{i,t}(z_{i,t}+\delta_tv_{i,t},x_{-i,t})-\nabla_if_{i,t}(z_t)\|]\nonumber\\
&\leq2B_xL\alpha_t\mathbf{E}_{v_{i,t}\in\mathbb{B}^{n_i}}[\|(\delta_tv_{i,t},x_{-i,t}-z_{-i,t})\|]\nonumber\\
&\leq2\sqrt{N}B_xL\alpha_t\delta_t,\label{equ55}
\end{align}
where $z_{-i,t}=col(z_{1,t},\ldots,z_{i-1,t},z_{i+1,t},\ldots,z_{N,t})$, the equality is based on the Leibniz integral rule, the first inequality is based on the Cauchy-Schwarz inequality and (\ref{equ7}), and the second inequality holds by Assumption \ref{assump2}.

For the third term on the right-hand side of (\ref{equ54}), one has that
\begin{align}
&\alpha_t\left\langle\breve{x}_{i,t}^*-{z}_{i,t},\nabla_if_{i,t}(\breve{x}_{t}^*)-{\nabla}_if_{i,t}({x}_{t}^*)\right\rangle\nonumber\\
&\leq\|\breve{x}_{i,t}^*-{z}_{i,t}\|\|\nabla_if_{i,t}(\breve{x}_{t}^*)-{\nabla}_if_{i,t}({x}_{t}^*)\|\nonumber\\
&\leq2B_xL\alpha_t\|\breve{x}_t^*-x_t^*\|\nonumber\\
&\leq2\sqrt{N}B_x^2L\alpha_t\eta_t.\label{equ56}
\end{align}

For the fourth term on the right-hand side of (\ref{equ54}), it is derived that
\begin{align}
&\alpha_t\left\langle\breve{x}_{i,t}^*-{z}_{i,t},{\nabla}_if_{i,t}(x_t^*)\right\rangle\nonumber\\
&=\alpha_t\left\langle{x}_{i,t}^*-{z}_{i,t},{\nabla}_if_{i,t}(x_t^*)+\gamma_t(\nabla{g}_{i,t}(x_{i,t}^*))^{\top}\lambda_t^*\right\rangle\nonumber\\
&~~~-\alpha_t\gamma_t\left\langle{x}_{i,t}^*-{z}_{i,t},(\nabla g_{i,t}(x_{i,t}^*))^{\top}\lambda_t^*\right\rangle\nonumber\\
&~~~-\alpha_t\eta_t\left\langle{x}_{i,t}^*,{\nabla}_if_{i,t}(x_t^*)\right\rangle\nonumber\\
&\leq2B_xL_g\Lambda\alpha_t\gamma_t+B_xL_f\alpha_t\eta_t,\label{equ57}
\end{align}
where the inequality is derived based on the optimality of $(x_{i,t}^*,\lambda_t^*)$ and (\ref{equ7}). Substituting (\ref{equ55}), (\ref{equ56}) and (\ref{equ57}) into (\ref{equ54}) implies that
\begin{align}
&\alpha_t\left\langle\breve{x}_{i,t}^*-{z}_{i,t},\nabla_i\hat{f}_{i,t}(z_t)\right\rangle\nonumber\\
&\leq\alpha_t\left\langle\breve{x}_{i,t}^*-{z}_{i,t},\nabla_if_{i,t}(z_t)-\nabla_if_{i,t}(\breve{x}_{t}^*)\right\rangle+2\sqrt{N}B_xL\alpha_t\delta_t\nonumber\\
&~~~+(2\sqrt{N}B_x^2L+B_xL_f)\alpha_t\eta_t+2B_xL_g\Lambda\alpha_t\gamma_t.\label{equ58}
\end{align}

For the second term on the right-hand side of (\ref{equ53}), there holds
\begin{align}
&\alpha_t\mathbf{E}_{u_{i,t}}\left[\left\langle{z}_{i,t}-\tilde{z}_{i,t+1},\hat{\nabla}_if_{i,t}(z_t)\right\rangle\right]\nonumber\\
&\leq\mathbf{E}_{u_{i,t}}\left[\frac{\alpha_t^2}{\mu_0}\|\hat{\nabla}_if_{i,t}(z_t)\|^2+\frac{\mu_0}{4}\|z_{i,t}-\tilde{z}_{i,t+1}\|^2\right]\nonumber\\
&\leq\frac{n_i^2B_f^2}{\mu_0}\frac{\alpha_t^2}{\delta_t^2}+\frac{\mu_0}{4}\mathbf{E}_{u_{i,t}}\left[\|z_{i,t}-\tilde{z}_{i,t+1}\|^2\right],\label{equ59}
\end{align}
where the first and second inequalities are derived by the Cauchy-Schwarz inequality and 5) in Lemma \ref{lemma1}, respectively.

As to the third term on the right-hand side of (\ref{equ53}), one gets that
\begin{align}
&\alpha_t\mathbf{E}_{u_{i,t}}\left[\left\langle\breve{x}_{i,t}^*-\tilde{z}_{i,t+1},(\hat{\nabla}g_{i,t}(z_{i,t}))^{\top}\tilde{\lambda}_{i,t}\right\rangle\right]\nonumber\\
&=\alpha_t\mathbf{E}_{u_{i,t}}\left[\left\langle{z}_{i,t}-\tilde{z}_{i,t+1},(\hat{\nabla}g_{i,t}(z_{i,t}))^{\top}\tilde{\lambda}_{i,t}\right\rangle\right]\nonumber\\
&~~~+\alpha_t\left[\left\langle\breve{x}_{i,t}^*-{z}_{i,t},({\nabla}\hat{g}_{i,t}(z_{i,t}))^{\top}\tilde{\lambda}_{i,t}\right\rangle\right]\nonumber\\
&\leq\mathbf{E}_{u_{i,t}}\left[\frac{1}{\mu_0}\alpha_t^2\|\tilde{\lambda}_{i,t}\|^2\|\hat{\nabla}g_{i,t}(z_{i,t})\|^2+\frac{\mu_0}{4}\|{z}_{i,t}-\tilde{z}_{i,t+1}\|^2\right]\nonumber\\
&~~~+\alpha_t\tilde{\lambda}_{i,t}^{\top}(\hat{g}_{i,t}(\breve{x}_{i,t}^*)-\hat{g}_{i,t}(z_{i,t}))\nonumber\\
&\leq\frac{n_i^2B_g^4}{\mu_0}\frac{\alpha_t^2}{\beta_t^2\delta_t^2}+\frac{\mu_0}{4}\mathbf{E}_{u_{i,t}}\left[\|{z}_{i,t}-\tilde{z}_{i,t+1}\|^2\right]\nonumber\\
&~~~+\alpha_t\overline{\lambda}_{t}^{\top}(\hat{g}_{i,t}(\breve{x}_{i,t}^*)-\hat{g}_{i,t}(z_{i,t}))\nonumber\\
&~~~+\alpha_t(\tilde{\lambda}_{i,t}-\overline{\lambda}_t)^{\top}(\hat{g}_{i,t}(\breve{x}_{i,t}^*)-\hat{g}_{i,t}(z_{i,t}))\nonumber
\end{align}
\begin{align}
&\leq\frac{n_i^2B_g^4}{\mu_0}\frac{\alpha_t^2}{\beta_t^2\delta_t^2}+\frac{\mu_0}{4}\mathbf{E}_{u_{i,t}}\left[\|{z}_{i,t}-\tilde{z}_{i,t+1}\|^2\right]\nonumber\\
&~~~+\alpha_t\overline{\lambda}_t^{\top}(g_{i,t}(x_{i,t}^*)-g_{i,t}(x_{i,t})+g_{i,t}(\breve{x}_{i,t}^*)-g_{i,t}(x_{i,t}^*))\nonumber\\
&~~~+\alpha_t\overline{\lambda}_t^{\top}(g_{i,t}(x_{i,t})-g_{i,t}(z_{i,t})+\hat{g}_{i,t}(\breve{x}_{i,t}^*)-g_{i,t}(\breve{x}_{i,t}^*))\nonumber\\
&~~~+\alpha_t\overline{\lambda}_t^{\top}(g_{i,t}(z_{i,t})-\hat{g}_{i,t}(z_{i,t}))+2B_g\alpha_t\|\tilde{\lambda}_{i,t}-\overline{\lambda}_t\|\nonumber\\
&\leq\frac{n_i^2B_g^4}{\mu_0}\frac{\alpha_t^2}{\beta_t^2\delta_t^2}+\frac{\mu_0}{4}\mathbf{E}_{u_{i,t}}\left[\|{z}_{i,t}-\tilde{z}_{i,t+1}\|^2\right]\nonumber\\
&~~~+\alpha_t\overline{\lambda}_t^{\top}(g_{i,t}(x_{i,t}^*)-g_{i,t}(x_{i,t}))+\alpha_t\|\overline{\lambda}_t\|\cdot\eta_tL_g\|x_{i,t}^*\|\nonumber\\
&~~~+\alpha_t\|\overline{\lambda}_t\|\cdot\delta_tL_g+\alpha_t\|\overline{\lambda}_t\|\cdot\delta_tL_g+2B_g\alpha_t\|\tilde{\lambda}_{i,t}-\overline{\lambda}_t\|\nonumber\\
&\leq\frac{n_i^2B_g^4}{\mu_0}\frac{\alpha_t^2}{\beta_t^2\delta_t^2}+\frac{\mu_0}{4}\mathbf{E}_{u_{i,t}}\left[\|{z}_{i,t}-\tilde{z}_{i,t+1}\|^2\right]\nonumber\\
&~~~+\alpha_t\overline{\lambda}_t^{\top}(g_{i,t}(x_{i,t}^*)-g_{i,t}(x_{i,t}))+B_xB_gL_g\frac{\alpha_t\eta_t}{\beta_t}\nonumber\\
&~~~+2B_gL_g\frac{\alpha_t\delta_t}{\beta_t}+4\sqrt{N}B_g^2\alpha_t\sum_{s=0}^{t-1}\sigma_m^s\gamma_{t-1-s},\label{equ60}
\end{align}
where the first equality is based on 2) of Lemma \ref{lemma1}, the first inequality is obtained by the Cauchy-Schwarz inequality and the convexity of $\hat{g}_{i,t}(\cdot)$, the other inequalities are based on Lemma \ref{lemma1}, and the last inequality is derived via Lemma \ref{lemma3}.


On the other hand, it can be derived from Assumptions \ref{assump5} and \ref{assump6} that
\begin{align}
&D_{\phi_i}(\breve{x}_{i,t+1}^*,z_{i,t+1})\nonumber\\
&\leq D_{\phi_i}(\breve{x}_{i,t}^*,z_{i,t+1})+K\|\breve{x}_{i,t+1}^*-\breve{x}_{i,t}^*\|,\label{equ61}\\
&D_{\phi_i}(\breve{x}_{i,t}^*,z_{i,t+1})\nonumber\\
&\leq(1-\alpha_t)D_{\phi_i}(\breve{x}_{i,t}^*,z_{i,t})+\alpha_tD_{\phi_i}(\breve{x}_{i,t}^*,\tilde{z}_{i,t+1}).\label{equ62}
\end{align}

Then, combining with (\ref{equ53}) and (\ref{equ58})--(\ref{equ62}) gives that
\begin{align}
&\mathbf{E}_{u_{i,t}}\left[D_{\phi_i}(\breve{x}_{i,t+1}^*,{z}_{i,t+1})\right]\nonumber\\
&\leq{D}_{\phi_i}(\breve{x}_{i,t}^*,z_{i,t})+K\|\breve{x}_{i,t+1}^*-\breve{x}_{i,t}^*\|\nonumber\\
&~~+\alpha_t^2\left\langle\breve{x}_{i,t}^*-{z}_{i,t},\nabla_if_{i,t}(z_t)-\nabla_if_{i,t}(\breve{x}_{t}^*)\right\rangle\nonumber\\
&~~+2\sqrt{N}B_xL\alpha_t^2\delta_t+2B_xL_g\Lambda\alpha_t^2\gamma_t\nonumber\\
&~~+(2\sqrt{N}B_x^2L+B_xL_f)\alpha_t^2\eta_t+\alpha_t^2\overline{\lambda}_t^{\top}(g_{i,t}(x_{i,t}^*)-g_{i,t}(x_{i,t}))\nonumber\\
&~~+\frac{n_i^2B_f^2}{\mu_0}\frac{\alpha_t^3}{\delta_t^2}+\frac{n_i^2B_g^4}{\mu_0}\frac{\alpha_t^3}{\beta_t^2\delta_t^2}+B_xB_gL_g\frac{\alpha_t^2\eta_t}{\beta_t}\nonumber\\
&~~+2B_gL_g\frac{\alpha_t^2\delta_t}{\beta_t}+4\sqrt{N}B_g^2\alpha_t^2\sum_{s=0}^{t-1}\sigma_m^s\gamma_{t-1-s}.\label{equ63}
\end{align}

Via the properties of conditional expectation, one has that
\begin{align*}
\mathbf{E}_{\mathcal{F}_T}\left[\mathbf{E}_{u_{i,t}}\left[D_{\phi_i}(\breve{x}_{i,t+1}^*,{z}_{i,t+1})\right]\right]
=\mathbf{E}\left[D_{\phi_i}(\breve{x}_{i,t+1}^*,{z}_{i,t+1})\right].
\end{align*}
For (\ref{equ63}), summing over $i=1,\ldots,N$, $t=1,\ldots,T$, taking expectation on its both sides, relying on Assumption \ref{assump3}, and rearranging the inequality, one has that
\begin{align}
&\mu\sum_{t=1}^T\mathbf{E}\left[\|\breve{x}_t^*-z_t\|^2\right]\nonumber\\
&\leq\sum_{t=1}^T\frac{1}{\alpha_t^2}\sum_{i=1}^N\mathbf{E}\left[D_{\phi_i}(\breve{x}_{i,t}^*,z_{i,t})-D_{\phi_i}(\breve{x}_{i,t+1}^*,z_{i,t+1})\right]\nonumber\\
&~~~+K\sum_{t=1}^T\frac{1}{\alpha_t^2}\sum_{i=1}^N\|\breve{x}_{i,t+1}^*-\breve{x}_{i,t}^*\|\nonumber\\
&~~~+2N\sqrt{N}B_xL\sum_{t=1}^T\delta_t+2NB_xL_g\Lambda\sum_{t=1}^T\gamma_t\nonumber\\
&~~~+N(2\sqrt{N}B_x^2L+B_xL_f)\sum_{t=1}^T\eta_t+\frac{n^2B_f^2}{\mu_0}\sum_{t=1}^T\frac{\alpha_t}{\delta_t^2}\nonumber\\
&~~~+\frac{n^2B_g^4}{\mu_0}\sum_{t=1}^T\frac{\alpha_t}{\beta_t^2\delta_t^2}-\sum_{t=1}^T\overline{\lambda}_t^{\top}g_{t}(x_{t})+NB_xB_gL_g\sum_{t=1}^T\frac{\eta_t}{\beta_t}\nonumber\\
&~~~+2NB_gL_g\sum_{t=1}^T\frac{\delta_t}{\beta_t}+4N\sqrt{N}B_g^2\sum_{t=1}^T\sum_{s=0}^{t-1}\sigma_m^s\gamma_{t-1-s}.\label{equ64}
\end{align}
where the inequality holds based on $\sum_{i=1}^Ng_{i,t}(x_{i,t}^*)\leq 0$.

Based on the definition of $\Lambda_t$ in Lemma \ref{lemma3}, the following equality holds:
\begin{align}
-\frac{\Lambda_{t+1}}{2\gamma_t}&=-\sum\limits_{i=1}^N\left[\frac{1}{2\gamma_t}\|\lambda_{i,t+1}-\lambda\|^2-\frac{1}{2\gamma_{t-1}}\|\lambda_{i,t}-\lambda\|^2\right]\nonumber\\
&~~~+\left(\frac{1}{2\gamma_t}-\frac{1}{2\gamma_{t-1}}-\frac{\beta_t}{2}\right)\sum_{i=1}^N\|\lambda_{i,t}-\lambda\|^2.
\end{align}
Summing over $t\in[T]$ yields that
\begin{align}
&-\sum_{t=1}^T\frac{\Delta_{t+1}}{2\gamma_t}\nonumber\\
&=-\sum\limits_{i=1}^N\left[\frac{1}{2\gamma_T}\|\lambda_{i,T+1}-\lambda\|^2-\frac{1}{2\gamma_{0}}\|\lambda_{i,1}-\lambda\|^2\right]\nonumber\\
&~~~+\sum\limits_{t=1}^T\left(\frac{1}{2\gamma_t}-\frac{1}{2\gamma_{t-1}}-\frac{\beta_t}{2}\right)\sum_{i=1}^N\|\lambda_{i,t}-\lambda\|^2\nonumber\\
&\leq\frac{N}{2}\|\lambda\|^2+\sum\limits_{t=1}^T\left(\frac{1}{2\gamma_t}-\frac{1}{2\gamma_{t-1}}-\frac{\beta_t}{2}\right)\sum_{i=1}^N\|\lambda_{i,t}-\lambda\|^2,\label{equ66}
\end{align}
where $\lambda_{i,1}={\bf0}_m$ is applied to derive the inequality.
Combining (\ref{equ40}) in Lemma \ref{lemma3}, (\ref{equ64}) and (\ref{equ66}) results in (\ref{equ48}). \hfill$\blacksquare$
\subsection{Proof of Lemma \ref{lemma5}}\label{C}
To derive the upper bounds of the expected dynamic regrets and the accumulated constraint violation, one first needs to analyze (\ref{equ48}). Note that
\begin{align}
&\sum_{t=1}^T\frac{1}{\alpha_t^2}\sum_{i=1}^N\left[D_{\phi_i}(\breve{x}_{i,t}^*,z_{i,t})-D_{\phi_i}(\breve{x}_{i,t+1}^*,z_{i,t+1})\right]\nonumber\\
&=\sum_{t=1}^T\sum_{i=1}^N\left[\frac{1}{\alpha_t^2}D_{\phi_i}(\breve{x}_{i,t}^*,z_{i,t})-\frac{1}{\alpha_{t+1}^2}D_{\phi_i}(\breve{x}_{i,t+1}^*,z_{i,t+1})\right]\nonumber\\
&~~~+\sum_{t=1}^T\sum_{i=1}^N\left[\frac{1}{\alpha_{t+1}^2}-\frac{1}{\alpha_t^2}\right]D_{\phi_i}(\breve{x}_{i,t+1}^*,z_{i,t+1})\nonumber\\
&\leq\frac{1}{\alpha_1^2}\sum_{i=1}^ND_{\phi_i}(\breve{x}_{i,1}^*,z_{i,1})+\sum_{t=1}^T\sum_{i=1}^N\left[\frac{1}{\alpha_{t+1}^2}-\frac{1}{\alpha_t^2}\right]2B_xK\nonumber\\
&\leq\frac{2NB_xK}{\alpha_{T+1}^2},
\end{align}
where (\ref{equ20}) is adopted to get the first inequality.

For the second term on the right-hand side of (\ref{equ48}), one has
\begin{align}
&K\sum_{t=1}^T\frac{1}{\alpha_t^2}\sum_{i=1}^N\|\breve{x}_{i,t+1}^*-\breve{x}_{i,t}^*\|\nonumber\\
&\leq\sum_{t=1}^T\frac{K}{\alpha_t^2}\sum_{i=1}^N(\|(1-\eta_{t+1})(x_{i,t+1}^*-x_{i,t}^*)+(\eta_t-\eta_{t+1})x_{i,t}^*\|)\nonumber\\
&\leq\frac{\sqrt{N}K}{\alpha_T^2}\Phi_T^*+NKB_x\sum_{t=1}^T\frac{\eta_t-\eta_{t+1}}{\alpha_t^2}.
\end{align}

Thus, by (\ref{equ48}), the upper bound of $\mu\sum_{t=1}^T\mathbf{E}\left[\|\breve{x}_t^*-z_t\|^2\right]$ can be obtained as
\begin{align}
&\mu\sum_{t=1}^T\mathbf{E}\left[\|\breve{x}_t^*-z_t\|^2\right]\nonumber\\
&\leq\frac{2NB_xK}{\alpha_{T+1}^2}+\frac{\sqrt{N}K}{\alpha_T^2}\Phi_T^*+NKB_x\sum_{t=1}^T\frac{\eta_t-\eta_{t+1}}{\alpha_t^2}\nonumber\\
&~~~+B_{1,T}+B_{\lambda,T}.\label{equ73}
\end{align}
Then, (\ref{equ67}) is derived by taking $\lambda=\mathbf{0}_m$ in (\ref{equ73}) and
\begin{align}
&\mathbf{E}\left[Reg_i(T)\right]\nonumber\\
&=\mathbf{E}\left[\sum_{t=1}^T(f_{i,t}(x_{i,t},x_{-i,t}^*)-f_{i,t}(x_{i,t}^*,x_{-i,t}^*))\right]\nonumber\\
&\leq{L}_f\mathbf{E}\left[\sum_{t=1}^T\|x_{i,t}-x_{i,t}^*\|\right]\nonumber\\
&={L}_f\mathbf{E}\left[\sum_{t=1}^T\|z_{i,t}-\breve{x}_{i,t}^*+\delta_tu_{i,t}-\eta_t\breve{x}_{i,t}^*\|\right]\nonumber\\
&\leq{L}_f\sum_{t=1}^T\mathbf{E}\left[\|z_{i,t}-\breve{x}_{i,t}^*\|\right]+L_f\sum_{t=1}^T(B_x\eta_t+\delta_t)\nonumber\\
&\leq L_f\sqrt{T\sum_{t=1}^T\mathbf{E}\left[\|z_{t}-\breve{x}_{t}^*\|^2\right]}+L_f\sum_{t=1}^T(B_x\eta_t+\delta_t).\label{equ72}
\end{align}

On the other hand, let $\lambda=\lambda_c=2B_{4,T}^{-1}[\sum_{t=1}^Tg_t(x_t)]_+$, then together with (\ref{equ73}), it follows that
\begin{align}
\frac{1}{B_{4,T}}\left\|\left[\sum_{t=1}^Tg_t(x_t)\right]_+\right\|^2&\leq{B}_{2,T}+\frac{\sqrt{N}K}{\alpha_T^2}\Phi_T^*+\Upsilon_{2}.
\end{align}
The proof is completed. \hfill$\blacksquare$


\subsection{Proof of Lemma \ref{lemma8}}\label{E}
For any $i\in[N]$, it can be derived based on the optimality of $\tilde{z}_{i,t+1}$ in (\ref{e54b}) that for any $z\in(1-\eta_t)X_i$,
\begin{align}
&\left\langle\tilde{z}_{i,t+1}-z,\alpha_t(\hat{\nabla}_if_{i,t-\tau}(z_{t-\tau})+(\hat{\nabla} g_{i,t-\tau}(z_{i,t-\tau}))^{\top}\tilde{\lambda}_{i,t})\right\rangle\nonumber\\
&+\left\langle\tilde{z}_{i,t+1}-z,\nabla\phi_i(\tilde{z}_{i,t+1})-\nabla\phi_i({z}_{i,t-\tau})\right\rangle\leq 0.\label{equ101}
\end{align}

Taking $z=\breve{x}^*_{i,t-\tau}$ in (\ref{equ101}) yields that
\begin{align}
&\alpha_t\langle\tilde{z}_{i,t+1}-\breve{x}_{i,t-\tau}^*,\hat{\nabla}_if_{i,t-\tau}(z_{t-\tau})+(\hat{\nabla}g_{i,t-\tau}(z_{i,t-\tau}))^{\top}\tilde{\lambda}_{i,t}\rangle\nonumber\\
&\leq\left\langle\breve{x}_{i,t-\tau}^*-\tilde{z}_{i,t+1},\nabla\phi_i(\tilde{z}_{i,t+1})-\nabla\phi_i({z}_{i,t-\tau})\right\rangle\nonumber\\
&=D_{\phi_i}(\breve{x}_{i,t-\tau}^*,z_{i,t-\tau})-D_{\phi_i}(\breve{x}_{i,t-\tau}^*,\tilde{z}_{i,t+1})\nonumber\\
&~~~-D_{\phi_i}(\tilde{z}_{i,t+1},z_{i,t-\tau})\nonumber\\
&\leq{D}_{\phi_i}(\breve{x}_{i,t-\tau}^*,z_{i,t-\tau})-D_{\phi_i}(\breve{x}_{i,t-\tau}^*,\tilde{z}_{i,t+1})\nonumber\\
&~~~-\frac{\mu_0}{2}\|\tilde{z}_{i,t+1}-z_{i,t-\tau}\|^2.\label{}
\end{align}
Then, one has
\begin{align}
&D_{\phi_i}(\breve{x}_{i,t-\tau}^*,\tilde{z}_{i,t+1})\nonumber\\
&\leq\alpha_t\left\langle\breve{x}_{i,t-\tau}^*-\tilde{z}_{i,t+1},\hat{\nabla}_if_{i,t-\tau}(z_{t-\tau})\right\rangle\nonumber\\
&~~~+\alpha_t\left\langle\breve{x}_{i,t-\tau}^*-\tilde{z}_{i,t+1},(\hat{\nabla}g_{i,t-\tau}(z_{i,t-\tau}))^{\top}\tilde{\lambda}_{i,t}\right\rangle\nonumber\\
&~~~+{D}_{\phi_i}(\breve{x}_{i,t-\tau}^*,z_{i,t-\tau})-\frac{\mu_0}{2}\|\tilde{z}_{i,t+1}-z_{i,t-\tau}\|^2.\label{}
\end{align}

On the other hand, it can be derived from Assumptions \ref{assump5} and \ref{assump6} that
\begin{align}
&D_{\phi_i}(\breve{x}_{i,t+1}^*,z_{i,t+1})\nonumber\\
&\leq D_{\phi_i}(\breve{x}_{i,t-\tau}^*,z_{i,t+1})+K\|\breve{x}_{i,t+1}^*-\breve{x}_{i,t-\tau}^*\|,\label{equ61}\\
&D_{\phi_i}(\breve{x}_{i,t-\tau}^*,z_{i,t+1})\nonumber\\
&\leq(1-\alpha_t)D_{\phi_i}(\breve{x}_{i,t-\tau}^*,z_{i,t-\tau})+\alpha_tD_{\phi_i}(\breve{x}_{i,t-\tau}^*,\tilde{z}_{i,t+1}).\label{equ62}
\end{align}
Then, following the proof of Lemma \ref{lemma4}, Lemma \ref{lemma8} can be proved.
\hfill$\blacksquare$

\subsection{Proof of Lemma \ref{lemma9}}\label{F}
Note that the first term on the right hand-side of (\ref{equ91}) satisfies
\begin{align}
&\sum_{t=\tau+1}^{T+\tau}\frac{1}{\alpha_t^2}\sum_{i=1}^N\left[D_{\phi_i}(\breve{x}_{i,t-\tau}^*,z_{i,t-\tau})-D_{\phi_i}(\breve{x}_{i,t+1}^*,z_{i,t+1})\right]\nonumber\\
&=\sum_{t=\tau+1}^{T+\tau}\sum_{i=1}^N\left[\frac{D_{\phi_i}(\breve{x}_{i,t-\tau}^*,z_{i,t-\tau})}{\alpha_{t-\tau-1}^2}-\frac{D_{\phi_i}(\breve{x}_{i,t+1}^*,z_{i,t+1})}{\alpha_{t}^2}\right]\nonumber\\
&~~~+\sum_{t=\tau+1}^{T+\tau}\sum_{i=1}^N\left[\frac{1}{\alpha_{t}^2}-\frac{1}{\alpha_{t-\tau-1}^2}\right]D_{\phi_i}(\breve{x}_{i,t-\tau}^*,z_{i,t-\tau})\nonumber\\
&\leq\sum_{\kappa=1}^{\tau+1}\sum_{i=1}^N\frac{D_{\phi_i}(\breve{x}_{i,\kappa}^*,z_{i,\kappa})}{\alpha_{\kappa-1}^2}+\sum_{t=\tau+1}^{T+\tau}\left[\frac{2NB_xK}{\alpha_{t}^2}-\frac{2NB_xK}{\alpha_{t-\tau-1}^2}\right]\nonumber\\
&\leq2NB_xK\sum_{\kappa=1}^{\tau+1}\frac{1}{\alpha_{T+\kappa-1}^2}\nonumber\\
&\leq\frac{2NB_xK(\tau+1)}{\alpha_{T+\tau}^2},\label{equ105}
\end{align}
where Assumption \ref{assump5} is adopted to get the first inequality. In addition, it holds that
\begin{align}
&\sum_{t=\tau+1}^{T+\tau}\frac{K}{\alpha_t^2}\sum_{i=1}^N\|\breve{x}_{i,t+1}^*-\breve{x}_{i,t-\tau}^*\|\nonumber\\
&\leq\sum_{t=\tau+1}^{T+\tau}\frac{K}{\alpha_t^2}\sum_{\kappa=0}^{\tau}\sum_{i=1}^N\|\breve{x}_{i,t+1-\kappa}^*-\breve{x}_{i,t-\kappa}^*\|\nonumber\\
&\leq\sum_{t=\tau+1}^{T+\tau}\frac{K}{\alpha_t^2}\sum_{\kappa=0}^{\tau}\sum_{i=1}^N\|(1-\eta_{t+1-\kappa})({x}_{i,t+1-\kappa}^*-x^*_{i,t-\kappa})\|\nonumber\\
&~~~+\sum_{t=\tau+1}^{T+\tau}\frac{K}{\alpha_t^2}\sum_{\kappa=0}^{\tau}\sum_{i=1}^N\|(\eta_{t-\kappa}-\eta_{t+1-\kappa}){x}_{i,t-\kappa}^*\|\nonumber\\
&\leq\frac{\sqrt{N}K}{\alpha_{T+\tau}^2}\sum_{t=\tau+1}^{T+\tau}\sum_{\kappa=0}^{\tau}\|{x}_{t+1-\kappa}^*-x^*_{t-\kappa}\|\nonumber\\
&~~~+NB_xK\sum_{t=\tau+1}^{T+\tau}\sum_{\kappa=0}^{\tau}\frac{\eta_{t-\kappa}-\eta_{t+1-\kappa}}{\alpha_t^2}\nonumber\\
&\leq\frac{\sqrt{N}K}{\alpha_{T+\tau}^2}(\tau+1)\Phi_{T+\tau}^*+NB_xK\sum_{t=\tau+1}^{T+\tau}\sum_{\kappa=0}^{\tau}\frac{\eta_{t-\tau}-\eta_{t+1-\tau}}{\alpha_t^2},\label{equ106}
\end{align}
where the last inequality is based on the decrease of $\eta_t-\eta_{t+1}$ with $t$.
Combining (\ref{equ91}), (\ref{equ105}) and (\ref{equ106}), it can be obtained from (\ref{equ72}) that (\ref{equ92}) is satisfied by letting $\lambda={\bf 0}_m$ in (\ref{equ91}).

On the other hand, taking $\lambda=\lambda_\tau=2C_{4,T}^{-1}\left[\sum_{t=1}^Tg_t(x_t)\right]_+$ in (\ref{equ91}), it holds that
\begin{align*}
\frac{1}{C_{4,T}}\left\|\left[\sum_{t=1}^Tg_t(x_t)\right]_+\right\|^2&\leq{C}_{2,T}+\frac{\sqrt{N}K}{\alpha_{T+\tau}^2}(\tau+1)\Phi_{T+\tau}^*\\
&~~~+\Upsilon_{\tau,2}.
\end{align*}
The proof is thus completed. \hfill$\blacksquare$
\end{appendix}
\bibliographystyle{IEEEtran}

\begin{thebibliography}{10}
\providecommand{\url}[1]{#1}
\csname url@samestyle\endcsname
\providecommand{\newblock}{\relax}
\providecommand{\bibinfo}[2]{#2}
\providecommand{\BIBentrySTDinterwordspacing}{\spaceskip=0pt\relax}
\providecommand{\BIBentryALTinterwordstretchfactor}{4}
\providecommand{\BIBentryALTinterwordspacing}{\spaceskip=\fontdimen2\font plus
\BIBentryALTinterwordstretchfactor\fontdimen3\font minus
  \fontdimen4\font\relax}
\providecommand{\BIBforeignlanguage}[2]{{%
\expandafter\ifx\csname l@#1\endcsname\relax
\typeout{** WARNING: IEEEtran.bst: No hyphenation pattern has been}%
\typeout{** loaded for the language `#1'. Using the pattern for}%
\typeout{** the default language instead.}%
\else
\language=\csname l@#1\endcsname
\fi
#2}}
\providecommand{\BIBdecl}{\relax}
\BIBdecl

\bibitem{ghaderi2014opinion}
J.~Ghaderi and R.~Srikant, ``{Opinion dynamics in social networks with stubborn
  agents: Equilibrium and convergence rate},'' \emph{Automatica}, vol.~50,
  no.~12, pp. 3209--3215, 2014.

\bibitem{saad2012game}
W.~Saad, Z.~Han, H.~V. Poor, and T.~Ba\c{s}ar, ``Game-theoretic methods for the
  smart grid: An overview of microgrid systems, demand-side management, and
  smart grid communications,'' \emph{IEEE Signal Processing Magazine}, vol.~29,
  no.~5, pp. 86--105, 2012.

\bibitem{stankovic2012distributed}
M.~S. Stankovic, K.~H. Johansson, and D.~M. Stipanovic, ``{Distributed seeking
  of Nash equilibria with applications to mobile sensor networks},'' \emph{IEEE
  Transactions on Automatic Control}, vol.~57, no.~4, pp. 904--919, 2012.

\bibitem{basar1999dynamic}
T.~Ba\c{s}ar and G.~J. Olsder, \emph{{Dynamic Noncooperative Game
  Theory}}.\hskip 1em plus 0.5em minus 0.4em\relax SIAM, 1999, vol.~23.

\bibitem{facchinei2010generalized}
F.~Facchinei and C.~Kanzow, ``{Generalized Nash equilibrium problems},''
  \emph{Annals of Operations Research}, vol. 175, no.~1, pp. 177--211, 2010.

\bibitem{shamma2005dynamic}
J.~S. Shamma and G.~Arslan, ``{Dynamic fictitious play, dynamic gradient play,
  and distributed convergence to Nash equilibria},'' \emph{IEEE Transactions on
  Automatic Control}, vol.~50, no.~3, pp. 312--327, 2005.

\bibitem{de2019distributed}
C.~De~Persis and S.~Grammatico, ``{Distributed averaging integral Nash
  equilibrium seeking on networks},'' \emph{Automatica}, vol. 110, p. 108548,
  2019.

\bibitem{gadjov2019passivity}
D.~Gadjov and L.~Pavel, ``{A passivity-based approach to Nash equilibrium
  seeking over networks},'' \emph{IEEE Transactions on Automatic Control},
  vol.~64, no.~3, pp. 1077--1092, 2019.

\bibitem{koshal2016distributed}
J.~Koshal, A.~Nedi{\'c}, and U.~V. Shanbhag, ``Distributed algorithms for
  aggregative games on graphs,'' \emph{Operations Research}, vol.~64, no.~3,
  pp. 680--704, 2016.

\bibitem{salehisadaghiani2019distributed}
F.~Salehisadaghiani, W.~Shi, and L.~Pavel, ``{Distributed Nash equilibrium
  seeking under partial-decision information via the alternating direction
  method of multipliers},'' \emph{Automatica}, vol. 103, pp. 27--35, 2019.

\bibitem{ye2021distributed}
M.~Ye, ``Distributed nash equilibrium seeking for games in systems with bounded
  control inputs,'' \emph{IEEE Transactions on Automatic Control}, vol.~66,
  no.~8, pp. 3833--3839, 2021.

\bibitem{lu2020online}
K.~Lu, H.~Li, and L.~Wang, ``{Online distributed algorithms for seeking
  generalized Nash equilibria in dynamic environments},'' \emph{IEEE
  Transactions on Automatic Control}, vol.~66, no.~5, pp. 2289--2296, 2021.

\bibitem{meng2021decentralized}
M.~Meng, X.~Li, Y.~Hong, J.~Chen, and L.~Wang, ``Decentralized online learning
  for noncooperative games in dynamic environments,'' \emph{arXiv preprint
  arXiv:2105.06200}, 2021.

\bibitem{bravo2018bandit}
M.~Bravo, D.~Leslie, and P.~Mertikopoulos, ``Bandit learning in concave
  $n$-person games,'' in \emph{Proceedings of the 32nd International Conference
  on Neural Information Processing Systems}, 2018, pp. 5666--5676.

\bibitem{mertikopoulos2019learning}
P.~Mertikopoulos and Z.~Zhou, ``Learning in games with continuous action sets
  and unknown payoff functions,'' \emph{Mathematical Programming}, vol. 173,
  no.~1, pp. 465--507, 2019.

\bibitem{zinkevich2003online}
M.~Zinkevich, ``Online convex programming and generalized infinitesimal
  gradient ascent,'' in \emph{Proceedings of the 20th International Conference
  on Machine Learning}, 2003, pp. 928--936.

\bibitem{hall2015online}
E.~C. Hall and R.~M. Willett, ``Online convex optimization in dynamic
  environments,'' \emph{IEEE Journal of Selected Topics in Signal Processing},
  vol.~9, no.~4, pp. 647--662, 2015.

\bibitem{salehisadaghiani2016distributed}
F.~Salehisadaghiani and L.~Pavel, ``{Distributed Nash equilibrium seeking: A
  gossip-based algorithm},'' \emph{Automatica}, vol.~72, pp. 209--216, 2016.

\bibitem{facchinei2009nash}
F.~Facchinei and J.-S. Pang, ``Nash equilibria: The variational approach,''
  \emph{Convex Optimization in Signal Processing and Communications}, pp.
  443--449, 2010.

\bibitem{pavel2020distributed}
L.~Pavel, ``{Distributed GNE seeking under partial-decision information over
  networks via a doubly-augmented operator splitting approach},'' \emph{IEEE
  Transactions on Automatic Control}, vol.~65, no.~4, pp. 1584--1597, 2020.

\bibitem{li2020}
X.~Li, L.~Xie, and Y.~Hong, ``Distributed aggregative optimization over
  multi-agent networks,'' \emph{IEEE Transactions on Automatic Control}, 2021,
  DoI: 10.1109/TAC.2021.3095456.

\bibitem{bauschke2001joint}
H.~H. Bauschke and J.~M. Borwein, ``{Joint and separate convexity of the
  Bregman distance},'' \emph{Studies in Computational Mathematics}, vol.~8, pp.
  23--36, 2001.

\bibitem{flaxman2005online}
A.~D. Flaxman, A.~T. Kalai, and H.~B. McMahan, ``Online convex optimization in
  the bandit setting: Gradient descent without a gradient,'' in
  \emph{Proceedings of the 16th Annual ACM-SIAM Symposium on Discrete
  Algorithms}, 2005, pp. 385--394.

\bibitem{yi2020distributed}
X.~Yi, X.~Li, L.~Xie, and K.~H. Johansson, ``Distributed online convex
  optimization with time-varying coupled inequality constraints,'' \emph{IEEE
  Transactions on Signal Processing}, vol.~68, pp. 731--746, 2020.

\bibitem{nedic2009approximate}
A.~Nedi{\'c} and A.~Ozdaglar, ``Approximate primal solutions and rate analysis
  for dual subgradient methods,'' \emph{SIAM Journal on Optimization}, vol.~19,
  no.~4, pp. 1757--1780, 2009.

\bibitem{yi2020distributed_b}
X.~Yi, X.~Li, T.~Yang, L.~Xie, T.~Chai, and K.~H. Johansson, ``Distributed
  bandit online convex optimization with time-varying coupled inequality
  constraints,'' \emph{IEEE Transactions on Automatic Control}, vol.~66,
  no.~10, pp. 4620--4635, 2021.

\bibitem{cao2021decentralized}
X.~Cao and T.~Basar, ``Decentralized online convex optimization with feedback
  delays,'' \emph{IEEE Transactions on Automatic Control}, 2021, DOI:
  10.1109/TAC.2021.3092562.

\end{thebibliography}


\begin{IEEEbiography}[{\includegraphics[width=1in,height=1.25in,clip,keepaspectratio]{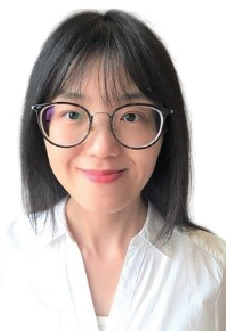}}]{Min Meng}
received the B.S. and Ph.D. degrees from Shandong University, China, in 2010 and 2015, respectively. She had a position as a Research Associate in Department of Mechanical Engineering, The University of Hong Kong, Hong Kong, China, from April to October in 2014, from July to September in 2016, and from January to March in 2017. From July 2015 to June 2016, she was a Research Associate in Department of Biomedical Engineering, City University of Hong Kong, Hong Kong, China. From July 2017 to September 2020, she worked as a Research Fellow in the School of Electrical and Electronic Engineering, Nanyang Technological University, Singapore. In 2020, she joined Tongji University, Shanghai, China, where she is now a professor.

Her research interests include multi-agent systems, distributed games and optimization, Boolean networks, distributed secure control and estimation, etc.
\end{IEEEbiography}

\begin{IEEEbiography}[{\includegraphics[width=1in,height=1.25in,clip,keepaspectratio]{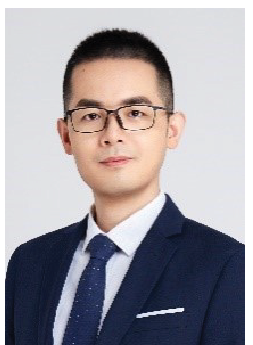}}]{Xiuxian Li}
(SM'21) received the B.S. degree in mathematics and applied mathematics and the M.S. degree in pure mathematics from Shandong University, Jinan, China, in 2009 and 2012, respectively, and the Ph.D. degree in mechanical engineering from the University of Hong Kong, Hong Kong, in 2016. From 2016 to 2020, he has been a research fellow with the School of Electrical and Electronic Engineering, Nanyang Technological University, Singapore, and he has also been a senior research associate with the Department of Biomedical Engineering, City University of Hong Kong, Hong Kong, in 2018. He held a visiting position at King Abdullah University of Science and Technology, Saudi Arabia, in September 2019. In 2020, he joined Tongji University, Shanghai, China, where he is now a professor.

His research interests include distributed control and optimization, algorithms, game theory, and machine learning, with applications to UAVs and autonomous vehicles, etc.
\end{IEEEbiography}

%

\begin{IEEEbiography}[{\includegraphics[width=1in,height=1.25in,clip,keepaspectratio]{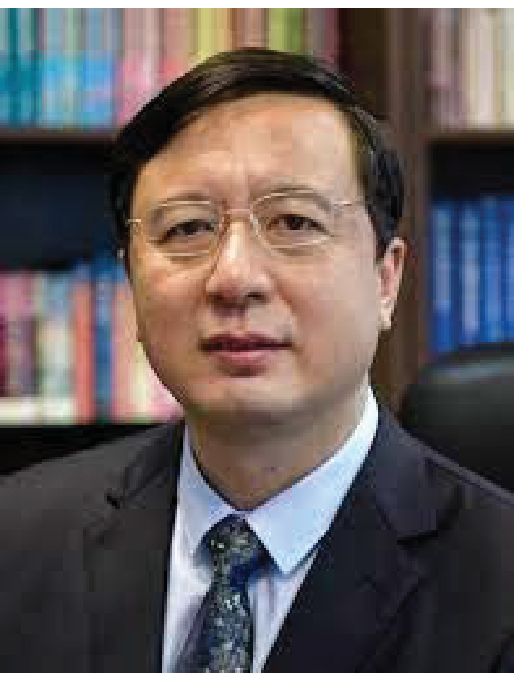}}]{Jie Chen} (F'19) received the B.S., M.S., and Ph.D. degrees in control theory and control engineering from the Beijing Institute of Technology (BIT), Beijing, China, in 1986, 1996, and 2001, respectively.
Dr. Chen is currently a Professor with the Department of Control Science and Engineering, Tongji University, Shanghai, China. He was a Professor with the School of Automation, BIT, and serves as the Director of the Key Laboratory of Intelligent Control and Decision of Complex Systems, BIT. His research interests include complex systems, multi-agent systems, multi-objective optimization and decision, constrained nonlinear control, and optimization methods.

Dr. Chen serves as the Vice Presidents of the Chinese Association of Automation (CAA) and the Chinese Association for Artificial Intelligence (CAAI). He serves as the Editor-in-Chief for Unmanned Systems and Autonomous Intelligent Systems, the Managing Editor for the Journal of Systems Science and Complexity, and an Editorial Board Member and an Associate Editor for several journals, including IEEE Transactions on Cybernetics, International Journal of Robust and Nonlinear Control, and Science China Information Sciences. He is an Academician of the Chinese Academy of Engineering, Fellow of IEEE, IFAC, CAA, and CAAI.
\end{IEEEbiography}
\end{document}